\documentclass[12pt]{article}
\usepackage{amssymb}
\usepackage{amsmath}
\usepackage[doublespacing]{setspace}
\usepackage{geometry}
\usepackage[final,dvips]{graphicx}

\newtheorem{theorem}{Theorem}

\newtheorem{corollary}{Corollary}

\newtheorem{lemma}{Lemma}

\newtheorem{proposition}{Proposition}

\input{tcilatex}
\begin{document}

\title{Detection of weak signals in high-dimensional complex-valued data}
\author{ Alexei Onatski\thanks{%
University of Cambridge, ao319@cam.ac.uk. }}
\maketitle

\begin{abstract}
This paper considers the problem of detecting a few signals in
high-dimensional complex-valued Gaussian data satisfying Johnstone's (2001) 
\textit{spiked covariance model}. We focus on the difficult case where
signals are weak in the sense that the sizes of the corresponding covariance
spikes are below the \textit{phase transition threshold} studied in Baik et
al (2005). We derive a simple analytical expression for the maximal possible
asymptotic probability of correct detection holding the asymptotic
probability of false detection fixed. To accomplish this derivation, we
establish what we believe to be a new formula for the \textit{%
Harish-Chandra/Itzykson-Zuber (HCIZ) integral} $\int_{\mathcal{U}\left(
p\right) }e^{\limfunc{tr}\left( AUB\bar{U}^{\prime }\right) }\left(
dU\right) $, where $A$ has a deficient rank $r<p$. The formula links the
HCIZ integral over $\mathcal{U}\left( p\right) $ to an HCIZ integral over a
potentially much smaller unitary group $\mathcal{U}\left( r\right) $. We
show that the formula generalizes to the integrals over orthogonal and
symplectic groups. In the most general form, it expresses the hypergeometric
function $_{0}F_{0}^{(\alpha )}$of two $p\times p$ matrix arguments as a
repeated contour integral of the hypergeometric function $_{0}F_{0}^{(\alpha
)}$of two $r\times r$ matrix arguments.
\end{abstract}

\textsc{Key words}: spiked covariance, sub-critical regime, signal
detection, sphericity tests, asymptotic power, contiguity, power envelope,
Harish-Chandra/Itzykson-Zuber integral, torus scalar product, hypergeometric
function.\bigskip

\section{Introduction}

Much contemporary research in statistics concerns with situations where the
dimensionality of data is large and comparable to the number of observations
(see special issues of the \textit{Philosophical Transactions of the Royal
Society} (2009) \textbf{367} and \textit{Annals of Statistics} (2008) 
\textbf{36}). Often, the goal is to estimate or detect a few signals
contaminated by high-dimensional noise. One general conclusion that seems to
emerge from this research is that, in the absence of \textit{a priori}
sparsity assumptions about signals, there is a lower limit for the
signal-to-noise ratio below which statistical inference about the signals
completely fails (Johnstone and Titterington, 2009, Nadakuditi and Edelman,
2008, Nadakuditi and Silverstein, 2010). This limit equals the \textit{phase
transition threshold} studied in Baik et al (2005). In a recent paper,
Onatski et al (2012) show that not all is lost below the threshold. They
consider the case of a single non-sparse signal in high-dimensional noisy
data and establish sharp non-trivial limits for the asymptotic power, as
both the data dimensionality and the number of observations go to infinity,
of statistical tests for signal detection when the signal may be arbitrarily
weak.

This paper extends Onatski et al (2012) to the case of multiple non-sparse
arbitrarily weak signals when the data are complex-valued. Complex-valued
data are of interest in signal processing (Schreier and Scharf, 2010),
wireless communication (Telatar, 1999, Tulino and Verdu, 2004), and the
spectral analysis of economic and financial time series (Onatski, 2009).
Considering the case of multiple signals is important for applied work
because the constraint that there is no more than one signal can rarely be
justified in practice. We derive a simple analytical expression for the
maximal possible asymptotic probability of correct detection, based on the
sample covariance eigenvalues of the data, holding the asymptotic
probability of false detection fixed.

We find that the asymptotic probability of detection may be close to one
even in cases where the strength of all signals is substantially below the
phase transition threshold. This finding is, perhaps, surprising in light of
the fact (P\'{e}ch\'{e}, 2003) that in such cases, sometimes referred to as
the sub-critical regime, the asymptotic behavior of any finite number of the
largest sample covariance eigenvalues is not different from their behavior
when the data are pure noise. We show that in these difficult cases, the
detection power lies not in the different behavior of a few of the largest
eigenvalues, but in the small deviations of the empirical distribution of
all the eigenvalues from the \textit{Marchenko-Pastur} \textit{limit}
(Marchenko and Pastur, 1967).

Let us discuss our findings in more detail. We assume that data consist of $%
n $ independent observations of $p$-dimensional complex-valued Gaussian
vectors $X_{t}$ with mean zero and covariance matrix $\sigma ^{2}\left(
I_{p}+VH\bar{V}^{\prime }\right) $, where $I_{p}$ is the $p$-dimensional
identity matrix, $\sigma $ is a real scalar, $H$ is an $r\times r$ real
diagonal matrix with elements $h_{j}\geq 0$ along the diagonal, and $V$ is a 
$\left( p\times r\right) $-dimensional complex parameter normalized so that $%
\bar{V}^{\prime }V=I_{r}$. Such a \textit{spiked covariance\ model} was
proposed by Johnstone (2001) as a simple model of a situation, often
observed in applications, where a few eigenvalues of the sample covariance
matrix, corresponding to signals, are relatively large, whereas the rest of
the eigenvalues are relatively small and tightly clustered. In our notation,
the size of the spikes is regulated by the values of $h_{j},$ and the signal
space is spanned by the columns of matrix $V$.

Let $\lambda _{1}\geq \lambda _{2}\geq ...\geq \lambda _{p}$ be the ordered
eigenvalues of $XX^{\prime }/n,$ where $X=\left[ X_{1},...,X_{n}\right] $,
and let $\lambda =\left( \lambda _{1},...,\lambda _{m}\right) ,$ where $%
m=\min \left\{ n,p\right\} $. We are interested in the asymptotic power of
tests for signal detection based on the information contained in $\lambda $
when $p,n\rightarrow \infty $ so that $p/n\rightarrow c$ with $0<c<\infty $.
Our null hypothesis is $H_{0}:h_{1}=...=h_{r}=0$ (no signals), and our
alternative is $H_{1}:$ $h_{i}>0$ for some $i=1,...,r$. The matrix $V$ is
left as an unspecified nuisance parameter. In this framework, signal
detection tests can also be interpreted as tests of sphericity.

We consider both cases of specified and unspecified $\sigma ^{2}$. For the
purpose of brevity, in the introduction, we will discuss only the case of
specified $\sigma ^{2}=1$. First, we study the likelihood ratio $L\left(
h;\lambda \right) ,$ defined as the ratio of the densities of $\lambda $
corresponding to unrestricted $h$ and restricted $h=0$, the densities being
evaluated at the observed value of $\lambda $. We show that $L\left(
h;\lambda \right) $ can be represented in the form of the determinant of an $%
r\times r$ matrix with entries equal to contour integrals of elementary
functions. We use Laplace approximations to these contour integrals to show
that for any $\bar{h}$ such that $0<\bar{h}<\sqrt{c},$ with $\sqrt{c}$ being
the value of the phase transition threshold, the sequence of log-likelihood
processes $\{\ln L\left( h;\lambda \right) ;h\in \lbrack 0,\bar{h}]^{r}\}$
converges weakly to a Gaussian process\footnote{%
Here the index $\lambda $ in the notation $\mathcal{L}_{\lambda }(h)$ is
used to distinguish the limiting log-likelihood process in the case of
specified $\sigma ^{2}=1,$ from that in the case of unspecified $\sigma
^{2}, $ which we denote by $\mathcal{L}_{\mu }(h)$.} $\left\{ \mathcal{L}%
_{\lambda }(h);h\in \left[ 0,\bar{h}\right] ^{r}\right\} $ under the null
hypothesis as $n,p\rightarrow \infty $. The limiting process has mean $%
\mathrm{E}\left[ \mathcal{L}_{\lambda }\mathcal{(}h\mathcal{)}\right] =\frac{%
1}{2}\sum_{i,j=1}^{r}\ln \left( 1-h_{i}h_{j}/c\right) $ and autocovariance
function $\mathrm{Cov}\left( \mathcal{L}_{\lambda }\left( h\right) ,\mathcal{%
L}_{\lambda }\left( \tilde{h}\right) \right) =-\sum_{i,j=1}^{r}\ln \left(
1-h_{i}\tilde{h}_{j}/c\right) $. The established weak convergence of
statistical experiments implies, via Le Cam's first lemma (see van der Vaart
1998, p.88), that the joint distributions of the sample covariance
eigenvalues under the null and under alternatives associated with $h\in %
\left[ 0,\sqrt{c}\right) ^{r}$ are mutually contiguous.

An asymptotic power envelope for eigenvalue-based tests of $H_{0}$ against $%
H_{1}$ can be constructed using the Neyman-Pearson lemma and Le Cam's third
lemma. We show that, for tests of size $\alpha ,$ the maximum achievable
asymptotic power against a point alternative $h=(h_{1},...,h_{r})$ equals $%
1-\Phi \left[ \Phi ^{-1}\left( 1-\alpha \right) -\sqrt{W}\right] ,$ where $%
\Phi $ is the standard normal distribution function and $W=-\sum_{i,j=1}^{r}%
\ln \left( 1-h_{i}h_{j}/c\right) $. A preliminary analysis indicates that
the asymptotic power of the likelihood ratio test based on the information
contained in $\lambda $ is close to the asymptotic power envelope. In
contrast, we find that the asymptotic powers of various previously proposed
tests are well below the envelope.

The central technical result of this paper is the contour integral
representation of the likelihood ratio. To derive such a representation, we
establish what we believe to be a novel formula for the hypergeometric
functions of two matrix arguments $_{0}F_{0}^{(\alpha )}\left( A,B\right) ,$
where a $p\times p$ matrix $A$ has rank $r<p,$ so that, without loss of
generality, only its upper-left $r\times r$ block $\mathcal{A}$ is non-zero.
Such functions appear as a key term in the explicit expressions for the
joint density of the eigenvalues of Wishart matrices with spiked covariance
parameter. In Lemma 1, we show that%
\begin{equation}
\left. _{0}F_{0}^{(\alpha )}\left( A,B\right) \right. =\frac{1}{r!\left(
2\pi \mathrm{i}\right) ^{r}}\oint_{\mathcal{K}}\!...\!\oint_{\mathcal{K}%
}\left. _{0}F_{0}^{(\alpha )}\left( \mathcal{A},\mathcal{Z}\right) \right.
\omega ^{\left( \alpha \right) }\left( A,B,\mathcal{Z}\right)
\prod\limits_{j=1}^{r}\mathrm{d}z_{j},  \label{formula}
\end{equation}%
where $\mathcal{Z}=\mathrm{diag}\left( z_{1},...,z_{r}\right) $ is an
auxiliary matrix, and $\omega ^{\left( \alpha \right) }\left( \cdot \right) $
is a simple function of $A,B,$ and $\mathcal{Z}$. This formula expresses the
hypergeometric function of high-dimensional arguments as a repeated contour
integral of a hypergeometric function of low-dimensional arguments, which is
convenient for analysis.

For the special case $r=1,$ (\ref{formula}) reduces to the formula that has
been recently derived in Mo\ (2011) and, independently, in Wang (2012) and
Onatski et al (2012) (see also Forrester, 2011 for a short derivation). Our
method of proof is different from the methods used by these authors. It is
based on the orthogonality of Jack polynomials with respect to the torus
scalar product (Macdonald (1995), Chapter VI, \S 10).

Although our analysis of signal detection in complex data requires only the
formula for $_{0}F_{0}^{(1)}\left( A,B\right) ,$ we establish (\ref{formula}%
) for all $\alpha =2/\beta ,$ where $\beta $ is a positive integer.\footnote{%
For cases where $\beta $ is odd, we require that $p-r+1$ be even. For even $%
\beta ,$ such a requirement is not needed.} The importance of finding
\textquotedblleft serviceable approximations\textquotedblright\ to $%
_{0}F_{0}^{(2)}\left( A,B\right) $ has been recently emphasized by Johnstone
(2007, p.322). Since in applications that rely on the spiked covariance
matrix framework $r$ is typically much smaller than $p$, analyzing $%
_{0}F_{0}^{(2)}\left( \mathcal{A},\mathcal{Z}\right) $ is much easier than
analyzing $_{0}F_{0}^{(2)}\left( A,B\right) $, and the established contour
integral representation of the latter may be of the welcomed service to
practitioners.

For $\alpha =2,1\,$\ and $1/2,$ function $_{0}F_{0}^{(\alpha )}\left(
A,B\right) $ has an integral representation $\int_{\mathcal{G}^{(\alpha
)}\left( p\right) }e^{\limfunc{tr}\left( AGBG^{-1}\right) }\left( \mathrm{d}%
G\right) ,$ where $\mathcal{G}^{(\alpha )}\left( p\right) $ is the
orthogonal group $\mathcal{O}(p)$ for $\alpha =2$, the unitary group $%
\mathcal{U}(p)$ for $\alpha =1$, and the compact symplectic group $\mathcal{S%
}\mathit{p}\left( p\right) $ for $\alpha =1/2,$ and where $\left( \mathrm{d}%
G\right) $ is the normalized Haar measure over $\mathcal{G}^{(\alpha
)}\left( p\right) $. Such integrals have various important applications in
mathematics and physics, where they are referred to as
Harish-Chandra/Itzykson-Zuber (HCIZ) integrals (Zinn-Justin and Zuber,
2003). The HCIZ integrals with rank-deficient $A$ have been used in the
analysis of spin glasses (Marinari et al, 1994), wireless communication
systems (Muller et al, 2008), statistical tests for signal detection
(Bianchi et al, 2010, and Onatski et al, 2012), distribution of the largest
sample covariance eigenvalue (Mo, 2011, and Wang, 2012), and spiked Wishart $%
\beta $-ensembles (Forrester, 2011). Their asymptotic behavior as $%
p\rightarrow \infty $ has been studied in Guionnet and Ma\"{\i}da (2005) and
Collins and \'{S}niady (2007). We hope that the reduction of HCIZ integrals
over large group $\mathcal{G}^{(\alpha )}\left( p\right) $ to those over
smaller group $\mathcal{G}^{(\alpha )}\left( r\right) $ that follows from (%
\ref{formula}) will be useful in a wide spectrum of applications.

The rest of this paper is organized as follows. In Section 2, we derive
explicit formulae for the likelihood ratios. Section 3 establishes
relationship (\ref{formula}). Section 4 uses (\ref{formula}) to derive
contour integral representations for the likelihood ratios. Section 5
applies Laplace approximations to the contour integrals in the derived
representation to obtain the asymptotics of the likelihood ratio process.
This asymptotics is then used along with the Neyman-Pearson lemma and Le
Cam's third lemma to establish a simple analytical formula for the maximal
possible asymptotic probability of correct signal detection holding the
asymptotic probability of false detection fixed. Section 6 concludes. All
proofs are relegated to the Appendix.

\section{Likelihood ratios}

As mentioned above, we assume that data consist of $n$ independent
observations of $p$-dimensional complex-valued Gaussian vectors $X_{t}\sim
N_{\mathbb{C}}\left( 0,\Sigma \right) $. This means that $X_{t}=Y_{t}+%
\mathrm{i}Z_{t}$, where $\mathrm{i}$ denotes the imaginary unit, and the
joint density of $\left( Y_{t},Z_{t}\right) $ at $\left( y,z\right) $
equals\ $\frac{1}{\left( 2\pi \right) ^{p}\det \Sigma }\exp \left\{ -%
\limfunc{tr}\left[ \Sigma ^{-1}\left( y+\mathrm{i}z\right) \left( y-\mathrm{i%
}z\right) ^{\prime }\right] \right\} $ (see, for example, Goodman, 1963).
Further, we assume that the covariance matrix $\Sigma $ equals $\sigma
^{2}\left( I_{p}+VH\bar{V}^{\prime }\right) ,$ where $H=\mathrm{diag}\left(
h_{1},...,h_{r}\right) $ quantifies the sizes of the covariance spikes. Our
goal is to study the asymptotic power of tests of $H_{0}:h_{1}=...=h_{r}=0$
against $H_{1}:$ $h_{i}>0$ for some $i=1,...,r$.

If $\sigma ^{2}$ is specified, the model is invariant with respect to
unitary transformations and the maximal invariant statistic is $\lambda $,
the vector of the first $m=\min \left\{ n,p\right\} $ eigenvalues of $%
XX^{\prime }/n,$ where $X=\left[ X_{1},...,X_{n}\right] $. Therefore, we
consider tests based on $\lambda $. If $\sigma ^{2}$ is unspecified, the
model is invariant with respect to the unitary transformations and
multiplications by non-zero scalars, and the maximal invariant is the vector
of normalized eigenvalues $\mu =\left( \mu _{1},...,\mu _{m-1}\right) ,$
where $\mu _{j}=\lambda _{j}/\left( \lambda _{1}+...+\lambda _{p}\right) $.
Hence, we consider tests based on $\mu $. Note that the distribution of $\mu 
$ does not depend on $\sigma ^{2},$ whereas if $\sigma ^{2}$ is specified,
we can always normalize $\lambda $ dividing it by $\sigma ^{2}$. Therefore,
in what follows, we will assume without loss of generality that $\sigma
^{2}=1$.

Let $h=\left( h_{1},...,h_{r}\right) $, and let us denote the joint density
of $\lambda _{1},...,\lambda _{m}$ as $p_{\lambda }\left( x;h\right) ,$ $%
x=\left( x_{1},...,x_{m}\right) \in (\mathbb{R}^{+})^{m}$ and that of $\mu
_{1},...,\mu _{m-1}$ as $p_{\mu }\left( y;h\right) ,$ $y=\left(
y_{1},...,y_{m-1}\right) \in (\mathbb{R}^{+})^{m-1}.$ We have%
\begin{equation}
p_{\lambda }\left( x;h\right) =\tilde{\gamma}\frac{\prod_{i=1}^{m}x_{i}^{%
\left\vert p-n\right\vert }\prod_{i<j}^{m}\left( x_{i}-x_{j}\right) ^{2}}{%
\prod_{i=1}^{r}\left( 1+h_{i}\right) ^{n}}\int\limits_{\mathcal{U}\left(
p\right) }e^{-n\limfunc{tr}\left( \Pi G\mathcal{X}G^{-1}\right) }\left( 
\mathrm{d}G\right) ,  \label{common complex1}
\end{equation}%
where $\tilde{\gamma}$ depends only on $n$ and $p$; $\Pi =\mathrm{diag}%
\left( \left( 1+h_{1}\right) ^{-1},...,\left( 1+h_{r}\right)
^{-1},1,...,1\right) $; $\mathcal{X}=\mathrm{diag}\left(
x_{1},...,x_{m},0,...,0\right) $ is a $\left( p\times p\right) $ diagonal
matrix, so that there are no zeros along the diagonal if $m=p$; $\mathcal{U}%
\left( p\right) $ is the set of all $p\times p$ unitary matrices; and $%
\left( \mathrm{d}G\right) $ is the invariant measure on the unitary group $%
\mathcal{U}\left( p\right) $ normalized to make the total measure unity.
Formula (\ref{common complex1}) is a special case of the densities given in
James (1964, p.489) for $n\geq p$ and in Ratnarajah and Vaillancourt (2005)
for $n<p$.

Let $s=x_{1}+...+x_{m}$ and let $y_{j}=x_{j}/s.$ Note that the Jacobian of
the coordinate change from $\left( x_{1},...,x_{m}\right) $ to $\left(
y_{1},...,y_{m-1},s\right) $ equals $s^{m-1}.$ Changing variables in (\ref%
{common complex1}) and integrating $s$ out, we obtain%
\begin{equation}
p_{\mu }\left( y;h\right) =\tilde{\gamma}\frac{\prod_{i=1}^{m}y_{i}^{\left%
\vert p-n\right\vert }\prod_{i<j}^{m}\left( y_{i}-y_{j}\right) ^{2}}{%
\prod_{i=1}^{r}\left( 1+h_{i}\right) ^{n}}\int_{0}^{\infty
}s^{np-1}\int\limits_{\mathcal{U}\left( p\right) }e^{-ns\limfunc{tr}\left(
\Pi G\mathcal{Y}G^{-1}\right) }\left( \mathrm{d}G\right) \mathrm{d}s,
\label{common complex2}
\end{equation}%
where $\mathcal{Y}=\mathrm{diag}\left( y_{1},...,y_{m-1},0,...,0\right) $ is
a $\left( p\times p\right) $ diagonal matrix.

Consider the likelihood ratios: $L\left( h;\lambda \right) =p_{\lambda
}\left( \lambda ;h\right) /p_{\lambda }\left( \lambda ;0\right) $ and $%
L\left( h;\mu \right) =p_{\mu }\left( \mu ;h\right) /p_{\mu }\left( \mu
;0\right) $. Formulae (\ref{common complex1}) and (\ref{common complex2})
imply the following Proposition.

\begin{proposition}
\textit{Let }$\mathcal{U}\left( p\right) $ \textit{be the set of all }$%
p\times p$\textit{\ unitary matrices. Denote by }$\left( \mathrm{d}G\right) $%
\textit{\ the invariant measure on the unitary group} $\mathcal{U}\left(
p\right) $ \textit{normalized to make the total measure unity. Further, let }%
$\Lambda =\mathrm{diag}\left( \lambda _{1},...,\lambda _{p}\right) $ and $M=%
\mathrm{diag}\left( \mu _{1},...,\mu _{p}\right) $. \textit{Then}%
\begin{eqnarray}
L\left( h;\lambda \right) &=&\prod_{i=1}^{r}\left( 1+h_{i}\right)
^{-n}\int\limits_{\mathcal{U}\left( p\right) }e^{-n\limfunc{tr}\left( \left(
\Pi -I\right) G\Lambda G^{-1}\right) }\left( \mathrm{d}G\right) \text{ 
\textit{and}}  \label{LR1} \\
L\left( h;\mu \right) &=&\frac{\prod_{i=1}^{r}\left( 1+h_{i}\right)
^{-n}n^{np}}{\Gamma \left( np\right) }\int_{0}^{\infty
}\!\!s^{np-1}e^{-ns}\!\!\int\limits_{\mathcal{U}\left( p\right) }\!\!e^{-ns%
\limfunc{tr}\left( \left( \Pi -I\right) GMG^{-1}\right) }\left( \mathrm{d}%
G\right) \mathrm{d}s.  \label{LR2}
\end{eqnarray}
\end{proposition}

Our analysis of the asymptotic power of tests for signal detection is based
on a study of the asymptotic properties of the likelihood ratio processes $%
\left\{ L\left( h;\lambda \right) ;h\!\in \!\left( R^{+}\right) ^{r}\right\} 
$ and $\left\{ L\left( h;\mu \right) ;h\!\in \!\left( R^{+}\right)
^{r}\right\} $. First, we will focus on the key terms in the expressions (%
\ref{LR1}) and (\ref{LR2}), which are the integrals over the unitary group.
These integrals are special cases of the complex hypergeometric function $%
_{0}F_{0}^{(1)}\left( A,B\right) =\int_{\mathcal{U}\left( p\right) }e^{%
\limfunc{tr}\left( AGBG^{-1}\right) }\!\left( \mathrm{d}G\right) ,$ where $A$
and, possibly, $B$ are rank-deficient. In the next section, we derive a
formula for $_{0}F_{0}^{(\alpha )}\left( A,B\right) $ with rank-deficient $A$
and $B$ that links this function to a hypergeometric functions of full-rank
matrix arguments of lower dimensions. We do not restrict attention to the
case $\alpha =1$ because, as discussed in the introduction, other cases
constitute independent interest.

\section{Contour integral representation for $_{0}F_{0}^{(\protect\alpha %
)}\left( A,B\right) $}

Let us first provide a necessary background on hypergeometric functions. Let 
$A$ and $B$ be Hermitian $p\times p$ matrices over real, complex, or
quaternion division algebra. The eigenvalues of such matrices are real and
we will denote them as $a=\left( a_{1},...,a_{p}\right) $ and $b=\left(
b_{1},...,b_{p}\right) $. The hypergeometric function $_{0}F_{0}^{(\alpha
)}\left( A,B\right) $ is defined as (see, for example, Koev and Edelman,
2006)%
\begin{equation}
_{0}F_{0}^{(\alpha )}\left( A,B\right) =\sum_{k=0}^{\infty }\sum_{\kappa
\vdash k}\frac{1}{k!}\frac{C_{\kappa }^{\left( \alpha \right) }\left(
A\right) C_{\kappa }^{\left( \alpha \right) }\left( B\right) }{C_{\kappa
}^{\left( \alpha \right) }\left( I_{p}\right) },  \label{F00}
\end{equation}%
where $C_{\kappa }^{\left( \alpha \right) }\left( A\right) =C_{\kappa
}^{\left( \alpha \right) }\left( a\right) ,$ $C_{\kappa }^{\left( \alpha
\right) }\left( B\right) =C_{\kappa }^{\left( \alpha \right) }\left(
b\right) $ and $C_{\kappa }^{\left( \alpha \right) }\left( I_{p}\right)
=C_{\kappa }^{\left( \alpha \right) }\left( 1,...,1\right) $ are normalized
Jack polynomials (Macdonald, 1995, chapter VI, \S 10), and the inner sum
runs over all partitions $\kappa $ of $k,$ that is over all non-increasing
sequences of non-negative integers $\kappa =\left( \kappa _{1},\kappa
_{2},...\right) $ such that $\kappa _{1}+\kappa _{2}+...=k.$ The
normalization of $C_{\kappa }^{\left( \alpha \right) }\left(
x_{1},...,x_{p}\right) $ is chosen so that 
\begin{equation}
\left( x_{1}+...+x_{p}\right) ^{k}=\sum_{\kappa \vdash k}C_{\kappa }^{\left(
\alpha \right) }\left( x_{1},...,x_{p}\right) .  \label{normalization}
\end{equation}

Note that $_{0}F_{0}^{(\alpha )}\left( A,B\right) $ depends on $A$ and $B$
only through $a$ and $b.$ Therefore, in what follows, without loss of
generality, we will consider only diagonal matrices $A\!=\!\mathrm{diag}%
\left( a_{1},...,a_{p}\right) $ and $B\!=\!\mathrm{diag}\left(
b_{1},...,b_{p}\right) .$ We will allow $a_{j}$ and $b_{j}$ to be complex,
thus extending definition (\ref{F00}) to complex diagonal matrices $A$ and $%
B $.

As was mentioned in the introduction, for $\alpha =2,1$ and $1/2,$
hypergeometric functions $_{0}F_{0}^{(\alpha )}\left( A,B\right) $ admit the
integral representation%
\begin{equation}
_{0}F_{0}^{(\alpha )}\left( A,B\right) =\int_{\mathcal{G}^{(\alpha )}\left(
p\right) }e^{\limfunc{tr}\left( AGBG^{-1}\right) }\left( \mathrm{d}G\right) ,
\label{representation}
\end{equation}%
where $\mathcal{G}^{(\alpha )}\left( p\right) $ is the orthogonal group $%
\mathcal{O}\left( p\right) $ for $\alpha =2,$ the unitary group $\mathcal{U}%
\left( p\right) $ for $\alpha =1,$ and the compact symplectic group $%
\mathcal{S}\mathit{p}\left( p\right) $ for $\alpha =1/2$. For real diagonal $%
A$ and $B,$ such a representation follows from the fact that $\int_{\mathcal{%
G}^{(\alpha )}\left( p\right) }C_{\kappa }^{\left( \alpha \right) }\left(
AGBG^{-1}\right) \left( \mathrm{d}G\right) =\frac{C_{\kappa }^{\left( \alpha
\right) }\left( A\right) C_{\kappa }^{\left( \alpha \right) }\left( B\right) 
}{C_{\kappa }^{\left( \alpha \right) }\left( I_{p}\right) }$ (see
Proposition 5.5 of Gross and Richards, 1987), and the fact that $%
\sum_{k=0}^{\infty }\sum_{\kappa \vdash k}\frac{1}{k!}C_{\kappa }^{\left(
\alpha \right) }\left( AGBG^{-1}\right) =e^{tr\left( AGBG^{-1}\right) },$
which follows from (\ref{normalization}). For complex diagonal $A$ and $B$,
the representation holds by the analytic continuation because both parts of
equality (\ref{representation}) are complex analytic functions of the
diagonal elements of $A$ and $B$.

The main result of this section is as follows.

\begin{lemma}
Let $A\!=\!\mathrm{diag}\left( a_{1},...,a_{p}\right) $ and $B\!=\!\mathrm{%
diag}\left( b_{1},...,b_{p}\right) $, where $a_{j}$ and $b_{j}$ are real or
complex numbers. Assume that $a_{j}\neq 0$ for $1\leq j\leq r$ and $a_{j}=0$
for $r<j\leq p$, and denote the upper left block of $A,$ $\mathrm{diag}%
\left( a_{1},...,a_{r}\right) ,$ as $\mathcal{A}$. Further, let $\mathcal{Z=}%
\mathrm{diag}\left( z_{1},...,z_{r}\right) ,$ where $z_{j}$ are complex
variables, and let $\mathcal{K}$ be a contour in the complex plane that
encircles $b_{1},...,b_{p}$ counter-clockwise. Finally, let $\alpha =2/\beta
,$ where $\beta $ is a positive integer. Then, assuming that $p-r+1$ is an
even integer in cases where $\beta $ is odd, and without this additional
assumption in cases where $\beta $ is even, we have%
\begin{equation}
\left. _{0}F_{0}^{(\alpha )}\left( A,B\right) \right. =\frac{1}{r!\left(
2\pi \mathrm{i}\right) ^{r}}\oint_{\mathcal{K}}\!...\!\oint_{\mathcal{K}%
}\left. _{0}F_{0}^{(\alpha )}\left( \mathcal{A},\mathcal{Z}\right) \right.
\omega ^{\left( \alpha \right) }\left( A,B,\mathcal{Z}\right)
\prod\limits_{j=1}^{r}\mathrm{d}z_{j},  \tag{1}  \label{main formula}
\end{equation}%
where%
\begin{eqnarray*}
\omega ^{\left( \alpha \right) }\left( A,B,\mathcal{Z}\right) &=&\left(
-1\right) ^{r\left( r-1\right) /\left( 2\alpha \right) }\prod_{j=1}^{r}\left[
\frac{\Gamma \left( \left( p\!+\!1\!-\!j\right) \!/\!\alpha \right) \Gamma
\left( 1/\alpha \right) }{\Gamma \left( \left( r\!+\!1\!-\!j\right)
\!/\!\alpha \right) }\right] \times \\
&&\prod\limits_{j>i}^{r}\left( z_{j}\!-\!z_{i}\right) ^{2/\alpha
}\!\prod\limits_{j=1}^{r}\left[ a_{j}^{1\!-\!(p\!-\!r\!+\!1)/\alpha
}\!\prod_{s=1}^{p}\left( z_{j}\!-\!b_{s}\right) ^{-1/\alpha }\right]
\end{eqnarray*}
\end{lemma}

The proposition reduces $\left. _{0}F_{0}^{(\alpha )}\left( A,B\right)
\right. $, a hypergeometric function with potentially high-dimensional
matrix arguments, to a repeated contour integral of $\left.
_{0}F_{0}^{(\alpha )}\left( \mathcal{A},\mathcal{Z}\right) \right. $, a
hypergeometric function with matrix arguments of possibly much lower
dimensions. In the special case where $r=1,$ $\left. _{0}F_{0}^{(\alpha
)}\left( \mathcal{A},\mathcal{Z}\right) \right. =e^{a_{1}z_{1}}$ and (\ref%
{main formula}) becomes 
\begin{equation}
\left. _{0}F_{0}^{(\alpha )}\left( A,B\right) \right. =\Gamma \left(
p/\alpha \right) a_{1}^{1-p/\alpha }\frac{1}{2\pi \mathrm{i}}\oint_{\mathcal{%
K}}e^{a_{1}z_{1}}\prod_{s=1}^{p}\left( z_{1}\!-\!b_{s}\right) ^{-\frac{1}{%
\alpha }}\mathrm{d}z_{1}.  \label{simple formula}
\end{equation}

For $\alpha =2,$ this formula has been established by Mo (2011), who used it
to analyze the asymptotic behavior of the largest eigenvalue of a rank-one
perturbation of a real Wishart matrix. He gives two proofs of the formula.
One of the proofs uses Jack polynomial expansions and requires that $p$ be
an even integer (consistent with our requirement that $p-r+1$ is even). The
other proof, which Mo (2011) calls geometric, allows for odd $p$.

Similar to the first proof of Mo, our proof of Lemma 1 uses Jack polynomial
expansions. In contrast to that proof, we do not rely on the simplification
of the Jack polynomials for top-order partitions, but use Jack polynomials'
orthogonality with respect to the torus scalar product (Macdonald, chapter
VI, \S 10). It is likely that our requirement that $p-r+1$ is even in cases
where $\beta =2/\alpha $ is odd can be lifted without affecting relationship
(\ref{main formula}). This would require a different proof of the
proposition, which is left for future research.

For $\alpha =2$ and $\alpha =2/\beta $ with even $\beta $, formula (\ref%
{simple formula}) has been independently established by Wang (2012). He uses
the formula to study the asymptotic distribution of the largest eigenvalue
of the real, complex and quaternionic Wishart matrices perturbed by matrices
of rank one. Wang's proof is similar to the first proof of Mo (2011) (see
Forrester, 2011, for an alternative proof). For $\alpha =2,$ formula (\ref%
{simple formula}) has also been independently established by Onatski et al
(2012). Their proof is based on the properties of the so-called Lauricella
function.

In contrast to (\ref{simple formula}), the general relationship (\ref{main
formula}) contains special functions on both left- and right-hand sides.
However, for $\alpha =1,$ it is possible to further simplify the right-hand
side of (\ref{main formula}) using Harish-Chandra/Itzykson-Zuber formula
(see Harish-Chandra, 1957, and Itzykson and Zuber, 1980)%
\begin{equation}
\left. _{0}F_{0}^{(1)}\left( \mathcal{A},\mathcal{Z}\right) \right. =\frac{%
\prod_{j=1}^{r-1}j!}{V_{r}\left( \mathcal{A}\right) V_{r}\left( \mathcal{Z}%
\right) }\det_{1\leq i,j\leq r}\left( e^{a_{i}z_{j}}\right) \text{,}
\label{HCIZ-full1}
\end{equation}%
where $V_{r}\left( \mathcal{A}\right) \!=\!\prod_{j>i}^{r}\left(
a_{j}\!-\!a_{i}\right) $ and $V_{r}\left( \mathcal{Z}\right)
\!=\!\prod_{j>i}^{r}\left( z_{j}\!-\!z_{i}\right) $ are the Vandermonde
determinants associated with the diagonal elements $a_{1},...,a_{r}$ of $%
\mathcal{A}$ and the diagonal elements $z_{1},...,z_{r}$ of $\mathcal{Z}$,
respectively. Using (\ref{HCIZ-full1}) in (\ref{main formula}), noting that
one of the terms in the definition of $\omega ^{(1)}\left( A,B,\mathcal{Z}%
\right) $ equals $V_{r}\left( \mathcal{Z}\right) ^{2},$ and applying
Andreief's identity (Andreief, 1883)%
\begin{equation*}
\det_{1\leq i,j\leq r}\left( \int f_{i}\left( x\right) g_{j}(x)\mathrm{d}\mu
\left( x\right) \right) =\frac{1}{r!}\int ...\int \det_{1\leq i,j\leq
r}\left( f_{i}\left( x_{j}\right) \right) \det_{1\leq i,j\leq r}\left(
g_{i}\left( x_{j}\right) \right) \prod\limits_{j}\mathrm{d}\mu \left(
x_{j}\right) ,
\end{equation*}%
we obtain the following Corollary.

\begin{corollary}
Under assumptions of Lemma 1, 
\begin{equation}
\left. _{0}F_{0}^{(1)}\left( A,B\right) \right. \!=\!\frac{\!\left(
-1\right) ^{r\left( r-1\right) /2}\!}{V_{r}\left( \mathcal{A}\right) }%
\prod_{j=1}^{r}\frac{\!\left( p\!-\!j\right) !}{\!a_{j}^{p-r}}\det_{1\leq
i,j\leq r}\left( \!\frac{1}{2\pi \mathrm{i}}\oint_{\mathcal{K}}\!\frac{%
e^{a_{i}z}z^{j-1}\mathrm{d}z}{\prod_{s=1}^{p}\left( z\!-\!b_{s}\right) }%
\!\right) .  \label{complex formula}
\end{equation}
\end{corollary}

An alternative way of deriving (\ref{complex formula}) is to apply l'H\^{o}%
pital's rule to the Harish-Chandra/Itzykson-Zuber determinantal formula 
\begin{equation}
\left. _{0}F_{0}^{(1)}\left( A,B\right) \right. =\frac{\prod_{j=1}^{p-1}j!}{%
V_{p}\left( A\right) V_{p}\left( B\right) }\det_{1\leq i,j\leq p}\left(
e^{a_{i}b_{j}}\right) ,  \label{HCIZ-full2}
\end{equation}%
the right-hand side of which is degenerate because $A$ is rank-deficient. We
include a proof of (\ref{complex formula}) that uses this approach in the
Supplementary Appendix. The proof is elementary in the sense that it does
not rely on properties of Jack polynomials.

\section{Likelihood ratios as contour integrals}

Combining Proposition 1 and Corollary 1 leads to useful contour integral
representations of the likelihood ratios $L\left( h;\lambda \right) $ and $%
L\left( h;\mu \right) $. We now introduce new notation to express such
representations in a convenient form. For any $z\in \mathcal{K},$ let us
define a random variable%
\begin{equation}
\Delta _{p}\left( z\right) =\sum_{j=1}^{p}\ln \left( z-\lambda _{j}\right)
-p\int \ln \left( z-\lambda \right) \mathrm{d}\mathcal{F}_{p}\left( \lambda
\right) ,  \label{delta_definition}
\end{equation}%
where $\mathcal{F}_{p}\left( \lambda \right) $ is the cumulative
distribution function of the Marchenko-Pastur distribution with a point mass
of $\max \left( 0,1-c_{p}^{-1}\right) $ at zero, where $c_{p}=p/n,$ and
density%
\begin{equation}
\psi _{p}\left( x\right) =\frac{1}{2\pi c_{p}x}\sqrt{\left( \bar{b}%
_{p}-x\right) \left( x-\bar{a}_{p}\right) },  \label{Marchenko-Pastur}
\end{equation}%
where $\bar{a}_{p}=\left( 1-\sqrt{c_{p}}\right) ^{2}$ and $\bar{b}%
_{p}=\left( 1+\sqrt{c_{p}}\right) ^{2}$. Further, let 
\begin{eqnarray}
f_{i}(z) &=&-\left( \frac{h_{i}}{1+h_{i}}z-c_{p}\int \ln \left( z-\lambda
\right) \mathrm{d}\mathcal{F}_{p}\left( \lambda \right) \right) \text{ and}
\label{f definition} \\
g_{j}(z) &=&z^{j-1}\exp \left\{ -\Delta _{p}\left( z\right) \right\} .
\label{g definition}
\end{eqnarray}%
Finally, for any permutation $\rho $ of the sequence $\left(
1,2,...,r\right) $ and any vector $\mathbf{z}=\left( z_{1},...,z_{r}\right) $%
, let%
\begin{equation}
q_{\rho }\left( \mathbf{z}\right) =\left( 1-\sum_{j=1}^{r}\frac{h_{\rho (j)}%
}{1+h_{\rho (j)}}\frac{z_{j}}{S}\right) ^{-p\left( n-r\right) -r(r+1)/2}\exp
\left\{ -\sum_{j=1}^{r}\frac{nh_{\rho (j)}z_{j}}{1+h_{\rho (j)}}\right\} ,
\label{h_tauz}
\end{equation}%
where $S=\lambda _{1}+...+\lambda _{p}.$

\begin{proposition}
\textit{Let the contour }$\mathcal{K}$ that encircles $\lambda
_{1},...,\lambda _{p}$ counter-clockwise \textit{be chosen so that }for any $%
z\in \mathcal{K}$, $\func{Re}z<\left( \sum_{j=1}^{r}\frac{h_{j}}{1+h_{j}}%
\right) ^{-1}S$. Then%
\begin{eqnarray}
L\left( h;\lambda \right) &=&k_{1}\det_{1\leq i,j\leq r}\left( \frac{1}{2\pi 
\mathrm{i}}\oint_{\mathcal{K}}e^{-nf_{i}\left( z\right) }g_{j}\left(
z\right) \mathrm{d}z\right) \text{ and}  \label{LR contour c 1} \\
L\left( h;\mu \right) &=&k_{2}\sum_{\rho }\frac{\limfunc{sgn}\rho }{\left(
2\pi \mathrm{i}\right) ^{r}}\oint_{\mathcal{K}}...\oint_{\mathcal{K}}q_{\rho
}\left( \mathbf{z}\right) \prod_{j=1}^{r}\left\{ e^{-nf_{\rho \left(
j\right) }\left( z_{j}\right) }g_{j}\left( z_{j}\right) \right\} \mathrm{d}%
z_{r}...\mathrm{d}z_{1},  \label{LR contour c 2}
\end{eqnarray}%
where $\mathrm{i}$ denotes the imaginary unit, the summation in (\ref{LR
contour c 2}) is over all permutations $\rho $ of the sequence $\left(
1,2,...,r\right) $,%
\begin{eqnarray*}
k_{1} &=&\left( -1\right) ^{r\left( r-1\right) /2}n^{-pr+r\left( r+1\right)
/2}\prod_{i>j}^{r}\left( h_{i}-h_{j}\right) ^{-1}\prod_{t=1}^{r}\left[
h_{t}^{r-p}\left( 1+h_{t}\right) ^{p-n-1}\left( p-t\right) !\right] ,\text{
and} \\
k_{2} &=&k_{1}\left( nS\right) ^{pr-r(r+1)/2}\Gamma \left( p\left(
n-r\right) +r(r+1)/2\right) \left[ \Gamma \left( np\right) \right] ^{-1}.
\end{eqnarray*}
\end{proposition}

In the next section, we perform the asymptotic analysis of $L\left(
h;\lambda \right) $ and $L\left( h;\mu \right) $ that relies on the Laplace
approximations of the contour integrals in (\ref{LR contour c 1}) and~(\ref%
{LR contour c 2}) after the contours are suitably deformed without changing
the value of the integrals.

\section{Asymptotic analysis}

\begin{figure}[t!]
\centering
\includegraphics[width=4in]{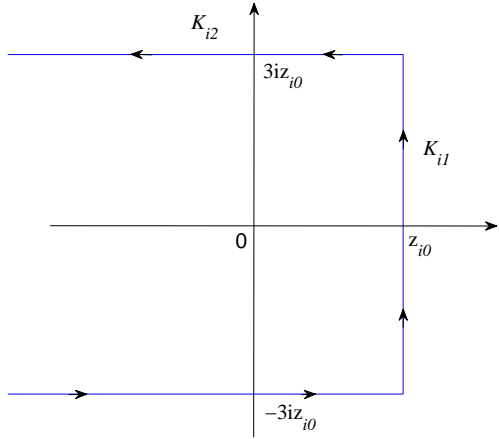}
\caption{Contour $\mathcal{K}_{i}$ (see (\protect\ref{contours})-(\protect\ref{contours 1})).}
\label{contour_distortionc}
\end{figure}

Consider contours $\mathcal{K}_{i}$ with $i=1,...,r$ which are obtained by
deforming the contour $\mathcal{K}$ defined in Proposition 2 so that $%
\mathcal{K}_{i}$ passes through 
\begin{equation*}
z_{i0}=\frac{\left( 1+h_{i}\right) \left( c_{p}+h_{i}\right) }{h_{i}}.
\end{equation*}%
Precisely, we define $\mathcal{K}_{i}$ as $\mathcal{K}_{i}=\mathcal{K}%
_{i+}\cup \mathcal{K}_{i-},$ where $\mathcal{K}_{i-}$ is the complex
conjugate of $\mathcal{K}_{i+}$ and $\mathcal{K}_{i+}=\mathcal{K}_{i1}\cup 
\mathcal{K}_{i2}$ with%
\begin{eqnarray}
\mathcal{K}_{i1} &\mathcal{=}&\left\{ z_{i0}+\mathrm{i}t:0\leq t\leq
3z_{i0}\right\} \text{ and}  \label{contours} \\
\text{ }\mathcal{K}_{i2} &\mathcal{=}&\left\{ x+3\mathrm{i}z_{i0}:-\infty
<x\leq z_{i0}\right\} .  \label{contours 1}
\end{eqnarray}%
Figure \ref{contour_distortionc} illustrates the choice of $\mathcal{K}_{i}$.

%

It is possible to verify that, when $0<h_{i}<\sqrt{c_{p}}$, the derivative
of $f_{i}\left( z\right) $ equals zero at $z_{i0}$. Therefore, choosing
contours of integration so they pass through $z_{i0}$ allows us to use the
method of steepest descent in the asymptotic analysis of the corresponding
integrals in (\ref{LR contour c 1}) and~(\ref{LR contour c 2}). The next
lemma shows that the change of contours in (\ref{LR contour c 1}) and (\ref%
{LR contour c 2}) does not lead to a change in the value of the
corresponding integrals.

\begin{lemma}
\textit{Suppose that the null hypothesis is true, and let }$\bar{h}$ \textit{%
be an arbitrary number such that }$0<\bar{h}<\sqrt{c}$. Suppose further that 
$h_{i}\leq \bar{h}$ for all $i=1,...,r$. Then, \textit{as }$n,p\rightarrow
\infty $ \textit{so that }$c_{p}\rightarrow c\in \left( 0,+\infty \right) ,$%
\begin{equation*}
\left( \oint_{\mathcal{K}}-\oint_{\mathcal{K}_{i}}\right) e^{-nf_{i}\left(
z\right) }g_{j}\left( z\right) \mathrm{d}z=0\text{ and}
\end{equation*}%
\begin{equation*}
\left( \oint_{\mathcal{K}}...\oint_{\mathcal{K}}-\oint_{\mathcal{K}_{\rho
(1)}}...\oint_{\mathcal{K}_{\rho (r)}}\right) q_{\rho }\left( \mathbf{z}%
\right) \prod_{j=1}^{r}\left\{ e^{-nf_{\rho \left( j\right) }\left(
z_{j}\right) }g_{j}\left( z_{j}\right) \right\} \mathrm{d}z_{r}...\mathrm{d}%
z_{1}=0.
\end{equation*}
\end{lemma}

Our next lemma establishes Laplace approximations to the contour integrals
in (\ref{LR contour c 1}) and (\ref{LR contour c 2}) after the change of the
contours. The lemma uses some new notation that we introduce now. When $%
f_{i}\left( z\right) $ is analytic at $z_{i0},$ let $f_{is}$ with $s=0,1,...$
be the coefficients in the power series representation%
\begin{equation}
f_{i}\left( z\right) =\sum_{s=0}^{\infty }f_{is}\left( z-z_{i0}\right) ^{s}.
\label{series for f and g}
\end{equation}%
When $f_{i}\left( z\right) $ is not analytic at $z_{i0}$, let the
coefficients $f_{is}$ be arbitrary numbers for all $s\in \mathbb{N}$.

\begin{lemma}
Under the conditions of Lemma 2, 
\begin{equation}
\oint\limits_{\mathcal{K}_{i}}e^{-nf_{i}(z)}g_{j}(z)\mathrm{d}z=e^{-nf_{i0}}%
\left[ \frac{g_{j}\left( z_{i0}\right) \pi ^{1/2}}{f_{i2}^{1/2}n^{1/2}}+%
\frac{O_{p}\left( 1\right) }{h_{i}^{j}n^{3/2}}\right] \text{ and}
\label{Watson}
\end{equation}%
\begin{eqnarray}
\oint_{\mathcal{K}_{\rho (1)}}... &&\!\!\!\oint_{\mathcal{K}_{\rho
(r)}}q_{\rho }\left( \mathbf{z}\right) \prod_{j=1}^{r}\left\{ e^{-nf_{\rho
\left( j\right) }\left( z_{j}\right) }g_{j}\left( z_{j}\right) \right\} 
\mathrm{d}z_{r}...\mathrm{d}z_{1}  \notag \\
&=&q_{\rho }\left( \mathbf{z}_{0}\right) \prod_{j=1}^{r}e^{-nf_{\rho \left(
j\right) 0}}\frac{g_{j}\left( z_{\rho \left( j\right) 0}\right) \pi ^{1/2}}{%
f_{\rho \left( j\right) 2}^{1/2}n^{1/2}}+\frac{O_{p}\left( 1\right) }{n}%
\prod_{j=1}^{r}\frac{e^{-nf_{\rho \left( j\right) 0}}}{h_{\rho \left(
j\right) }^{j}n^{1/2}},  \label{Watson1}
\end{eqnarray}%
where $O_{p}\left( 1\right) $\textit{\ is uniform in }$h_{1},...,h_{r}\in %
\left[ 0,\overline{h}\right] ,$ and $\mathbf{z}_{0}=\left( z_{\rho
(1)0},...,z_{\rho (r)0}\right) $.\textit{\ The branch of the square root in
formulae (\ref{Watson}) and (\ref{Watson1}) is chosen so that }$\left(
-1\right) ^{1/2}=-\mathrm{i}$.
\end{lemma}

Using Lemma 3, we establish the following theorem.

\begin{theorem}
\textit{Suppose that the null hypothesis is true (}$h=0$)\textit{. Let }$%
\bar{h}$\textit{\ be any fixed number such that }$0<\bar{h}<\sqrt{c}$ 
\textit{and let }$C\left[ 0,\overline{h}\right] ^{r}$\textit{\ be the space
of real-valued continuous functions on }$\left[ 0,\overline{h}\right] ^{r}$%
\textit{\ equipped with the supremum norm}. \textit{Then,} \textit{as }$%
n,p\rightarrow \infty $ \textit{so that }$p/n=c_{p}\rightarrow c\in \left(
0,+\infty \right) ,$\textit{\ we have}%
\begin{eqnarray}
L\!\left( h;\!\lambda \right) \!\!\!\! &=&\!\!\!\!\exp \!\left\{
\!-\!\sum\limits_{i=1}^{r}\!\Delta _{p}\!\left( \!z_{i0\!}\right) \!+\!\frac{%
1}{2}\!\sum\limits_{i,j=1}^{r}\ln \!\left( \!1\!-\!\frac{h_{i}h_{j}}{c_{p}}%
\!\right) \right\} \!+\!O_{p}\!\left( \!\frac{1}{n}\!\right) \text{ \textit{%
and}}  \label{equivalence 1} \\
L\!\left( h;\!\mu \right) \!\!\!\! &=&\!\!\!\!\exp \!\left\{
\!-\!\sum\limits_{i=1}^{r}\!\Delta _{p}\!\left( \!z_{i0}\!\right) \!+\!\frac{%
1}{2}\!\sum\limits_{i,j=1}^{r}\!\left( \!\ln \!\left( \!1\!\!-\!\!\frac{%
\!h_{i}h_{j}\!}{c_{p}}\!\right) \!+\!\frac{\!h_{i}h_{j}\!}{c_{p}}\!\right)
\!-\!\frac{\!S\!-\!p}{c_{p}}\!\sum\limits_{j=1}^{r}\!h_{j}\!\right\}
\!+\!O_{p}\!\!\left( \!\frac{1}{n}\!\right) \!,  \label{equivalence 2}
\end{eqnarray}%
where the $O_{p}\left( n^{-1}\right) $ terms are uniform in $h\in \left( 0,%
\bar{h}\right] ^{r}$. \textit{Furthermore, }$\ln L\left( h;\lambda \right) $ 
\textit{and} $\ln L\left( h;\mu \right) ,$\textit{\ viewed as random
elements of }$C\left[ 0,\overline{h}\right] ^{r}$,\textit{\ converge\ weakly
to }$\mathcal{L}_{\lambda }\left( h\right) $ \textit{and} $\mathcal{L}_{\mu
}\left( h\right) $ \textit{with Gaussian finite-dimensional distributions
such that, for any }$h,\tilde{h}\in \left[ 0,\overline{h}\right] ^{r},$%
\begin{eqnarray}
&&\limfunc{E}\left( \mathcal{L}_{\lambda }\left( h\right) \right) =\frac{1}{2%
}\sum\limits_{i,j=1}^{r}\ln \left( 1-\frac{h_{i}h_{j}}{c}\right) \text{,}
\label{mean} \\
&&\limfunc{Cov}\left( \mathcal{L}_{\lambda }\left( h\right) ,\mathcal{L}%
_{\lambda }\left( \tilde{h}\right) \right) =-\sum_{i,j=1}^{r}\ln \left( 1-%
\frac{h_{i}\tilde{h}_{j}}{c}\right) ,  \label{covariance} \\
&&\limfunc{E}\left( \mathcal{L}_{\mu }\left( h\right) \right) =\frac{1}{2}%
\sum_{i,j=1}^{r}\left( \ln \left( 1-\frac{h_{i}h_{j}}{c}\right) +\frac{%
h_{i}h_{j}}{c}\right) ,\text{ \textit{and}}  \label{mean mu} \\
&&\limfunc{Cov}\left( \mathcal{L}_{\mu }\left( h\right) ,\mathcal{L}_{\mu
}\left( \tilde{h}\right) \right) =-\sum_{i,j=1}^{r}\left( \ln \left( 1-\frac{%
h_{i}\tilde{h}_{j}}{c}\right) +\frac{h_{i}\tilde{h}_{j}}{c}\right) .
\label{covariance mu}
\end{eqnarray}
\end{theorem}

Theorem 1 and Le Cam's first lemma (van der Vaart (1998), p.88) imply that
the joint distributions of $\lambda _{1},...,\lambda _{m}$ (as well as those
of $\mu _{1},...,\mu _{m-1}$) under the null and under the alternative are
mutually contiguous for any $h\in \left[ 0,\sqrt{c}\right) ^{r}$. Along with
Le Cam's third lemma (van der Vaart (1998), p.90), this can be used to study
the \textquotedblleft local\textquotedblright\ powers of tests detecting
signals in noise.

Let $\beta _{\lambda }\left( h\right) $ and $\beta _{\mu }\left( h\right) $
be the asymptotic powers of the asymptotically most powerful $\lambda $- and 
$\mu $-based tests of size $\alpha $ of the null $h=0$ against a point
alternative $h=\left( h_{1},...,h_{r}\right) $ with $h_{j}<\sqrt{c},$ $%
j=1,...,r$. We have

\begin{theorem}
\textit{Let }$\Phi $\textit{\ denote the standard normal distribution
function. Then,}%
\begin{eqnarray}
\beta _{\lambda }\left( h\right) &=&1-\Phi \left[ \Phi ^{-1}\left( 1-\alpha
\right) -\sqrt{-\sum_{i,j=1}^{r}\ln \left( 1-\frac{h_{i}h_{j}}{c}\right) }%
\right] \text{ \textit{and}}  \label{local power} \\
\beta _{\mu }\left( h\right) &=&1-\Phi \left[ \Phi ^{-1}\left( 1-\alpha
\right) -\sqrt{-\sum_{i,j=1}^{r}\left( \ln \left( 1-\frac{h_{i}h_{j}}{c}%
\right) +\frac{h_{i}h_{j}}{c}\right) }\right] .  \label{local power mu}
\end{eqnarray}
\end{theorem}

The theorem implies in particular that detection of signals corresponding to
covariance spikes of sizes well below the phase transition threshold is
possible with high probability. Consider for example the case where the
number of observations equals the dimensionality of data so that $c=1,$ the
number of signals under the alternative equals five, and the signals have
equal but rather weak strengths $h_{1}=...=h_{5}=0.5$. Then the best
possible $\lambda $-based procedure for detecting such signals with the
asymptotic probability of false detection fixed at 0.05 has asymptotic
probability of correct detection $1-\Phi \left[ \Phi ^{-1}\left( 0.95\right)
-\sqrt{-25\ln \left( 1-0.25\right) }\right] \approx 0.85$.

Unfortunately, constructing testing procedures with uniformly optimal power
is hard because the log-likelihood process established in Theorem 1 is not
of the Gaussian shift type, so that the statistical experiments we study are
not locally asymptotically normal (LAN) ones. For the case of real-valued
data and $r=1,$ Onatski et al (2012) use numerical simulations to show that
the asymptotic powers of the likelihood ratio (LR) tests based on $\lambda $
and on $\mu $ are close to the respective asymptotic power envelopes $\beta
_{\lambda }\left( h\right) $ and $\beta _{\mu }\left( h\right) $. The $%
\lambda $- and $\mu $-based LR tests of $h=0$ against the alternative $h\in
\left( 0,\bar{h}\right) ^{r}$ reject the null if and only if , respectively, 
$2\sup_{h\in (0,\bar{h})}\ln L\left( h;\lambda \right) $ and $2\sup_{h\in (0,%
\bar{h})}\ln L\left( h;\mu \right) $ are sufficiently large. As $r$ grows,
it becomes increasingly difficult to find the asymptotic critical values for
the LR tests by simulation. This requires simulating an $r$-dimensional
Gaussian random field with the covariance function and the mean function
described in Theorem 1, which, for relatively large $r$, is computationally
expensive.

\begin{figure}[t!]
\centering
\includegraphics[width=5in]{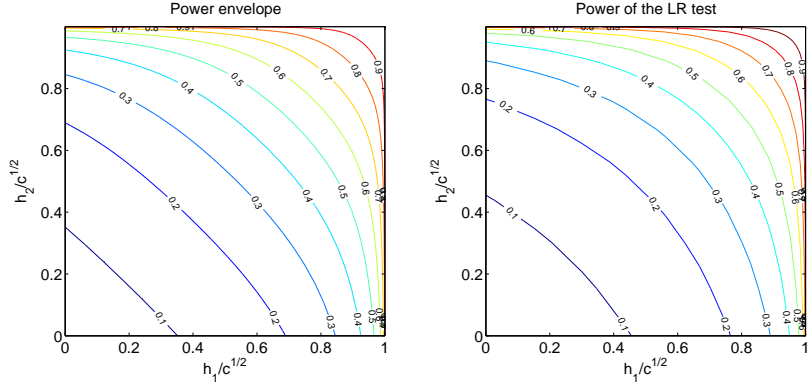}
\caption{The asymptotic power envelope $\protect\beta _{\protect\lambda }\left( h\right) $ and the asymptotic power
of the LR test based on $\protect\lambda $; $r=2,$ asymptotic size is $0.05$.}
\label{lrpower}
\end{figure}

For $r=2,$ Figure \ref{lrpower} shows the contour plots of the power
envelope $\beta _{\lambda }\left( h\right) $ (left panel) and of the
asymptotic power of the likelihood ratio test based on $\lambda $. We chose
parameter $\bar{h}$ so that it is very close to the threshold $\sqrt{c}$,
precisely $\bar{h}=\sqrt{c(1-e^{-36})}.$ We see that the contours of $\beta
_{\lambda }\left( h\right) $ and of the asymptotic power of the $\lambda $%
-based LR test corresponding to the same value of these functions are
relatively close to each other, which suggests that the LR test has good
asymptotic power properties. More detailed analysis of the asymptotic and
finite sample power of the LR test is, however, beyond the scope of this
paper, and is left for future research.

In contrast to the LR test, the popular signal detection procedures based on
the information in a few of the largest eigenvalues of $XX^{\prime }/n$
(see, for example, Krichman and Nadler (2009), Nadakuditi and Silverstein
(2010), Onatski (2009), Patterson et al (2006), Perry and Wolf (2010), and
Tracy and Widom (2009)), have trivial asymptotic power (that is, the
asymptotic power, which equals the asymptotic size) in the region $h\in %
\left[ 0,\sqrt{c}\right) ^{r}$. It is because the asymptotic behavior of any
finite number of the largest sample covariance eigenvalues when $h\in \left[
0,\sqrt{c}\right) ^{r}$ is not different from their behavior when the data
are pure noise (P\'{e}ch\'{e}, 2003).

As was mentioned above, signal detection tests can be interpreted as tests
of sphericity. Vice versa, previously proposed sphericity tests, can, in
principle, be used for signal detection. In the Supplementary Appendix, we
use Theorem 1 along with Le Cam's third lemma to derive asymptotic powers of
several such tests against \textquotedblleft spiked
covariance\textquotedblright\ alternatives. The derived asymptotic powers
turn out to be much lower than the asymptotic power envelopes $\beta
_{\lambda }\left( h\right) $ and $\beta _{\mu }\left( h\right) $. However,
we feel that this comparison is somewhat unfair to the sphericity tests
because they are typically designed against general alternatives, as opposed
to \textquotedblleft the spiked covariance\textquotedblright\ alternatives.
Therefore, and to save space, we do not report these results here.

\section{Conclusion}

This paper studies the asymptotic power of the signal detection tests in
complex-valued Gaussian data as both the number of observations and data
dimensionality go to infinity. Contrary to the conventional wisdom that
detection of signals becomes nearly impossible when their strength, measured
by the size of the covariance spikes, is below the phase transition
threshold, we find that detection of such signals may be possible with high
probability. The detection power lies not in the different behavior of a few
of the largest sample covariance eigenvalues under the null and the
alternative, which is exploited by the popular signal detection tests, but
in small deviations of the empirical distribution of all the eigenvalues
from the Marchenko-Pastur limit.

To derive our results, we consider the ratio of the densities of the sample
covariance eigenvalues under the null and under the alternative hypothesis.
We establish a contour integral representation of this likelihood ratio, and
use the Laplace approximation to derive its asymptotic limit. Our analysis
of the limiting log-likelihood ratio process shows that the sub-critical
region, where the sizes of the covariance spikes are below the phase
transition threshold, is the region of mutual contiguity of the joint
densities of the sample covariance eigenvalues under the null and the
alternative. We use the derived limiting log-likelihood process along with
Le Cam's third lemma and the Neyman-Pearson lemma to obtain the asymptotic
power envelope for the signal detection tests. Preliminary analysis
indicates that the asymptotic power of the likelihood ratio test based on
the sample covariance eigenvalues is close to the asymptotic power envelope.

Our technical analysis is based on what we believe to be a novel
representation of the Harish-Chandra/Itzykson-Zuber integral with one of the 
$p\times p$ matrices being of reduced rank $r$ in the form of an $r\times r$
matrix of contour integrals. We obtain such a representation as a corollary
to a much more general result established in Lemma 1. This result expresses
the hypergeometric function $_{0}F_{0}^{(\alpha )}$of two $p\times p$ matrix
arguments, one of which has rank $r$, as a repeated contour integral of the
hypergeometric function $_{0}F_{0}^{(\alpha )}$of two $r\times r$ matrix
arguments. As discussed in the introduction, the established dimension
reduction for the hypergeometric function may be important in various
applied and theoretical fields of study. In particular, for $\alpha =2$, it
can, potentially, be used to extend the analysis of this paper to the case
of real-valued data. Such an extension is currently under investigation.

\section{Appendix}

\textsc{Proof of Lemma 1.}

Let $f(\mathcal{Z)}$ and $g(\mathcal{Z})$ be functions defined on the $r$%
-dimensional torus $\left\{ \left\vert z_{j}\right\vert =1,\text{ for}%
\right. \newline
\left. \text{ }j=1,...,r\right\} .$ Consider the scalar product, sometimes
called the torus scalar product, 
\begin{equation}
\left\langle f,g\right\rangle _{\alpha }=\frac{1}{r!\left( 2\pi \mathrm{i}%
\right) ^{r}}\oint ...\oint f\left( \mathcal{Z}\right) \overline{g\left( 
\mathcal{Z}\right) }\prod_{i\neq j}\left( 1-z_{i}z_{j}^{-1}\right)
^{1/\alpha }\prod_{j=1}^{r}\frac{\mathrm{d}z_{j}}{z_{j}},  \label{torus}
\end{equation}%
where the contours of integration are the unit circles in the complex plane.
Our proof relies on the orthogonality property of Jack polynomials: $%
\left\langle C_{\kappa }^{(\alpha )},C_{\tau }^{(\alpha )}\right\rangle
_{\alpha }=0$\textit{\ for }$\kappa \neq \tau $ (Macdonald, chapter VI, \S %
10).

\begin{figure}[t!]
\centering
\includegraphics[width=4in]{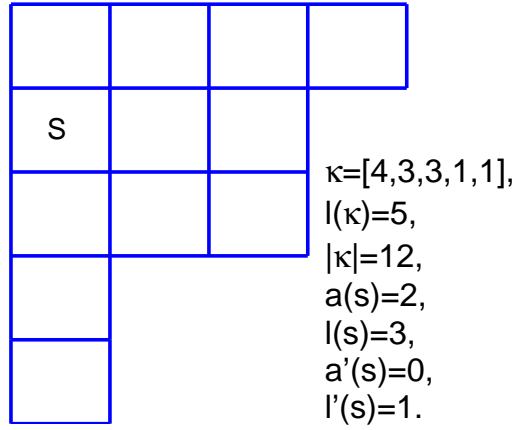}
\caption{The Ferrers diagram of partition 
$\left[ 4,3,3,1,1\right] $.}
\label{diagram}
\end{figure}%

Let us, first, introduce a few definitions (following Macdonald, 1995,
chapter I, \S 1, and Dumitriu et al, 2007): The non-zero $\kappa _{j}$ in
the partition $\kappa =\left[ \kappa _{1},\kappa _{2},...\right] $ are
called the parts of $\kappa $. The number of parts is the length of $\kappa
, $ denoted as $l\left( \kappa \right) .$ The sum of the parts is the weight
of $\kappa ,$ denoted as $\left\vert \kappa \right\vert $. We will identify
partition $\kappa $ with its Ferrers diagram, defined as an arrangement of $%
\left\vert \kappa \right\vert $ boxes in $l\left( \kappa \right) $
left-justified rows, the number of boxes in row $i$ being the same as $%
\kappa _{i}$ (see Figure \ref{diagram}). For each square $s\ $in the Ferrers
diagram, let $l^{\prime }\left( s\right) ,l(s),a(s),$ and $a^{\prime }(s)$
be respectively the numbers of squares in the diagram to the north, south,
east, and west of the square $s$. Further, let $h^{\ast }\left( s\right)
=l\left( s\right) +\alpha \left( 1+a(s)\right) $ and $h_{\ast }\left(
s\right) =l\left( s\right) +1+\alpha a(s).$ Finally, let $c\left( \kappa
,\alpha \right) =\prod_{s\in \kappa }h_{\ast }\left( s\right) $, $c^{\prime
}\left( \kappa ,\alpha \right) =\prod_{s\in \kappa }h^{\ast }\left( s\right) 
$, and $w\left( \kappa ,\alpha \right) =c\left( \kappa ,\alpha \right)
c^{\prime }\left( \kappa ,\alpha \right) $.

We will need the following lemmata.\medskip

\textbf{Lemma A1}. \textit{For the torus scalar product of }$C_{\kappa
}^{(\alpha )}$\textit{\ with itself, we have}%
\begin{equation}
\left\langle C_{\kappa }^{(\alpha )},C_{\kappa }^{(\alpha )}\right\rangle
_{\alpha }=\frac{\left( \alpha ^{\left\vert \kappa \right\vert }\left\vert
\kappa \right\vert !\right) ^{2}}{w\left( \kappa ,\alpha \right) }%
\prod_{j=1}^{r}\left( \frac{\Gamma \left( \frac{r-j+1}{\alpha }\right) }{%
\Gamma \left( \frac{1}{\alpha }\right) \Gamma \left( 1\!+\!\frac{r-j}{\alpha 
}\right) }\right) \prod_{s\in \kappa }\frac{r\!+\!a^{\prime }\left( s\right)
\alpha \!-\!l^{\prime }\left( s\right) }{r\!+\!\left( a^{\prime }\left(
s\right) \!+\!1\right) \alpha \!-\!l^{\prime }\left( s\right) \!-\!1}.
\label{scalar product}
\end{equation}%
\textsc{Proof:} Macdonald (1995, Chapter VI, \S 10) establishes the
orthogonality of \textquotedblleft $P$\textquotedblright\ normalizations of
Jack polynomials, $P_{\kappa }^{(\alpha )}$, with respect to the torus
scalar product. His formula (10.37) gives an explicit expression (up to a
constant that can be evaluated using (10.38)) for $\left\langle P_{\kappa
}^{(\alpha )},Q_{\kappa }^{(\alpha )}\right\rangle _{\alpha },$ where $%
Q_{\kappa }^{(\alpha )}=\frac{c\left( \kappa ,\alpha \right) }{c^{\prime
}\left( \kappa ,\alpha \right) }P_{\kappa }^{(\alpha )}$ (see (10.16)). On
the other hand, 
\begin{equation}
P_{\kappa }^{(\alpha )}=\frac{c^{\prime }\left( \kappa ,\alpha \right) }{%
\alpha ^{\left\vert \kappa \right\vert }\left\vert \kappa \right\vert !}%
C_{\kappa }^{(\alpha )}  \label{PClink}
\end{equation}%
(see, for example, Table 6 of Dumitriu et al, 2007). Substituting this
expression in Macdonald's formulae, we get (\ref{scalar product}).$\square $%
\medskip

\textbf{Lemma A2.} \textit{Let }$\mathcal{Z}=\mathrm{diag}\left(
z_{1},...,z_{r}\right) ,$\textit{\ where }$z_{1},...,z_{r}$\textit{\ are
complex variables, and let }$b_{1},...,b_{p}$\textit{\ be complex constants.
Then}%
\begin{equation}
\prod\limits_{j=1}^{r}\prod_{s=1}^{p}\left( 1\!-\!b_{s}z_{j}\right)
^{-1/\alpha }=\sum_{k=0}^{\infty }\sum_{\kappa \vdash k}\frac{w\left( \kappa
,\alpha \right) }{\left( \alpha ^{\left\vert \kappa \right\vert }\left\vert
\kappa \right\vert !\right) ^{2}}C_{\kappa }^{(\alpha )}\left( B\right)
C_{\kappa }^{(\alpha )}\left( \mathcal{Z}\right) .  \label{generating fn}
\end{equation}%
\textit{The series on the right-hand side of this equality converges
uniformly over }$\Omega _{\rho }=\left\{ \mathcal{Z}:\max_{j\leq
r}\left\vert z_{j}\right\vert \leq \rho ^{-1}\right\} ,$\textit{\ for any }$%
\rho >\max_{s\leq p}\left\vert b_{s}\right\vert $\textit{.}

\textsc{Proof:} Macdonald (1995, Chapter VI, \S 10) shows that $%
\prod\limits_{j=1}^{r}\prod\limits_{s=1}^{p}\left( 1\!-\!b_{s}z_{j}\right)
^{-1/\alpha }=\sum_{k=0}^{\infty }\sum_{\kappa \vdash k}P_{\kappa }^{(\alpha
)}\left( B\right) Q_{\kappa }^{(\alpha )}\left( \mathcal{Z}\right) ,$ where $%
Q_{\kappa }^{(\alpha )}=\frac{c\left( \kappa ,\alpha \right) }{c^{\prime
}\left( \kappa ,\alpha \right) }P_{\kappa }^{(\alpha )}$. This result
together with (\ref{PClink}) imply (\ref{generating fn}). The uniform
convergence in (\ref{generating fn}) follows from the fact that function $%
\prod_{j=1}^{r}\prod_{s=1}^{p}\left( 1\!-\!b_{s}z_{j}\right) ^{-1/\alpha }$
is analytic in an open region that includes $\Omega _{\rho }$.$\square
\medskip $

We are now ready to prove Lemma 1. Consider the right-hand side of (\ref%
{main formula}), which we will denote as RHS. We will assume that $%
\max_{s\leq p}\left\vert b_{s}\right\vert <1$ and that the contour $\mathcal{%
K}$ is the unit circle in the complex plane. That these assumptions are
without loss of generality follows from the fact that the value of RHS does
not change under the transformation $\mathcal{Z}\rightarrow \varphi \mathcal{%
Z},$ $B\rightarrow \varphi B,$ and $A\rightarrow \varphi ^{-1}A,$ where $%
\varphi $ is any positive number, and under a deformation of $\mathcal{K}$
into the unit circle (because such a deformation leaves the contour in the
region of the analyticity of the integrand). With these assumptions, and
noting that the component $\prod_{j>i}^{r}\left( z_{j}\!-\!z_{i}\right)
^{2/\alpha }$ of $\omega ^{(\alpha )}\left( A,B,\mathcal{Z}\right) $ equals $%
\left( -1\right) ^{r(r-1)/\left( 2\alpha \right)
}\prod_{j=1}^{r}z_{j}^{\left( r-1\right) /\alpha }\prod_{j\neq
i}^{r}\!\left( 1\!-\!z_{i}z_{j}^{-1}\right) ^{1/\alpha },$ we can rewrite
RHS for $\alpha =2/\beta ,$ where $\beta $ is any positive integer, in the
form of the torus scalar product%
\begin{equation*}
RHS=\gamma ^{(\alpha )}\left\langle \left. _{0}F_{0}^{(\alpha )}\!\left( \!%
\mathcal{A},\mathcal{Z}\!\right) \right. ,\prod\limits_{j=1}^{r}\left(
z_{j}/a_{j}\right) ^{(p\!-\!r\!+\!1)/\alpha
\!-\!1}\prod\limits_{j=1}^{r}\prod_{s=1}^{p}\left( 1\!-\!b_{s}z_{j}\right)
^{-1/\alpha }\right\rangle _{\alpha },
\end{equation*}%
where 
\begin{equation}
\gamma ^{(\alpha )}=\prod_{j=1}^{r}\left[ \frac{\Gamma \left( \left(
p\!+\!1\!-\!j\right) \!/\!\alpha \right) \Gamma \left( 1/\alpha \right) }{%
\Gamma \left( \left( r\!+\!1\!-\!j\right) \!/\!\alpha \right) }\right] .
\label{gamma_alfa}
\end{equation}

Substituting $_{0}F_{0}^{(\alpha )}\!\left( \!\mathcal{A},\mathcal{Z}%
\!\right) $ and $\prod\limits_{j=1}^{r}\prod_{s=1}^{p}\left(
1\!-\!b_{s}z_{j}\right) ^{-1/\alpha }$ in the above formula by their
expansions (\ref{F00}) and (\ref{generating fn}) in the series of Jack
polynomials, and interchanging the order of integration and summation, which
is possible because the series converge uniformly over the unit torus, we
obtain%
\begin{eqnarray*}
RHS &=&\gamma ^{(\alpha )}\sum_{k=0}^{\infty }\sum_{\kappa \vdash
k}\sum_{t=0}^{\infty }\sum_{\tau \vdash t}\frac{w\left( \tau ,\alpha \right) 
}{k!\left( \alpha ^{t}t!\right) ^{2}}\frac{C_{\kappa }^{\left( \alpha
\right) }\left( \mathcal{A}\right) C_{\tau }^{(\alpha )}\left( B\right) }{%
C_{\kappa }^{\left( \alpha \right) }\left( I_{r}\right) }\times \\
&&\left\langle C_{\kappa }^{\left( \alpha \right) }\left( \mathcal{Z}\right)
,\prod\limits_{j=1}^{r}\left( z_{j}/a_{j}\right) ^{(p\!-\!r\!+\!1)/\alpha
\!-\!1}C_{\tau }^{(\alpha )}\left( \mathcal{Z}\right) \right\rangle _{\alpha
}.
\end{eqnarray*}

But $\prod_{j=1}^{r}\left( z_{j}\right) ^{(p\!-\!r\!+\!1)/\alpha
\!-\!1}C_{\tau }^{(\alpha )}\left( \mathcal{Z}\right) =C_{\tilde{\tau}%
}^{(\alpha )}\left( \mathcal{Z}\right) ,$ where $\tilde{\tau}$ denotes
partition \linebreak $\lbrack \tau _{1}+\frac{p\!-\!r\!+\!1\!-\!\alpha }{%
\alpha },...,\tau _{r}+\frac{p\!-\!r\!+\!1\!-\!\alpha }{\alpha }]$. Note
that $\tilde{\tau}$ is well defined for $\alpha =2/\beta ,$ where $\beta $
is an even integer. If $\beta $ is an odd integer, we need to assume that $%
p\!-\!r\!+\!1$ is even. Therefore, using the orthogonality of the Jack
polynomials with respect to the torus scalar product, we have%
\begin{equation*}
RHS=\gamma ^{(\alpha )}\prod\limits_{j=1}^{r}a_{j}^{-(p\!-\!r\!+\!1)/\alpha
\!+\!1}\sum_{t=0}^{\infty }\sum_{\tau \vdash t}\frac{w\left( \tau ,\alpha
\right) }{\left\vert \tilde{\tau}\right\vert !\left( \alpha ^{t}t!\right)
^{2}}\frac{C_{\tilde{\tau}}^{\left( \alpha \right) }\!\left( \mathcal{A}%
\right) C_{\tau }^{(\alpha )}\!\left( B\right) }{C_{\tilde{\tau}}^{\left(
\alpha \right) }\!\left( I_{r}\right) }\left\langle C_{\tilde{\tau}}^{\left(
\alpha \right) }\!\left( \mathcal{Z}\right) ,\!C_{\tilde{\tau}}^{(\alpha
)}\!\left( \mathcal{Z}\right) \right\rangle _{\alpha }.
\end{equation*}%
Using Lemma A1, (\ref{gamma_alfa}), and equality $\prod%
\limits_{j=1}^{r}a_{j}^{-(p\!-\!r\!+\!1)/\alpha \!+\!1}C_{\tilde{\tau}%
}^{\left( \alpha \right) }\left( \mathcal{A}\right) =C_{\tau }^{\left(
\alpha \right) }\left( \mathcal{A}\right) =C_{\tau }^{\left( \alpha \right)
}\left( A\right) ,$ we get after some cancellations%
\begin{equation}
RHS=\sum_{t=0}^{\infty }\sum_{\tau \vdash t}\tilde{\gamma}^{(\alpha )}\frac{1%
}{t!}\frac{C_{\tau }^{\left( \alpha \right) }\left( A\right) C_{\tau
}^{(\alpha )}\left( B\right) }{C_{\tau }^{\left( \alpha \right) }\left(
I_{p}\right) },  \label{semi-final}
\end{equation}%
where%
\begin{equation*}
\tilde{\gamma}^{(\alpha )}=\frac{\alpha ^{2\left\vert \tilde{\tau}%
\right\vert }\left\vert \tilde{\tau}\right\vert !}{w\left( \tilde{\tau}%
,\alpha \right) }\frac{w\left( \tau ,\alpha \right) }{\alpha ^{2t}t!}\frac{%
C_{\tau }^{\left( \alpha \right) }\left( I_{p}\right) }{C_{\tilde{\tau}%
}^{\left( \alpha \right) }\left( I_{r}\right) }\prod_{j=1}^{r}\frac{\Gamma
\left( \left( p\!+\!1\!-\!j\right) \!/\!\alpha \right) }{\Gamma \left(
1\!+\!\left( r\!-\!j\right) \!/\!\alpha \right) }\prod_{s\in \tilde{\tau}}%
\frac{r\!+\!a^{\prime }\left( s\right) \alpha \!-\!l^{\prime }\left(
s\right) }{r\!+\!\left( a^{\prime }\left( s\right) \!+\!1\right) \alpha
\!-\!l^{\prime }\left( s\right) \!-\!1}.
\end{equation*}

In the above expression for $\tilde{\gamma}^{(\alpha )},$ substitute $%
C_{\tau }^{\left( \alpha \right) }\left( I_{p}\right) $ and $C_{\tilde{\tau}%
}^{\left( \alpha \right) }\left( I_{r}\right) $ by their explicit forms,
that can be obtained from a general formula%
\begin{equation}
C_{\kappa }^{\left( \alpha \right) }\left( I_{m}\right) =\frac{\alpha
^{\left\vert \kappa \right\vert }\left\vert \kappa \right\vert !}{w\left(
\kappa ,\alpha \right) }\prod_{s\in \kappa }\left( m+\alpha a^{\prime
}\left( s\right) -l^{\prime }(s)\right) .  \label{C of Identity}
\end{equation}%
A variant of this formula, that uses the generalized Pochhammer symbol, can
be found, for example, in Dumitriu et al (2007, Table 5). Then, after
cancellations, we get%
\begin{equation*}
\tilde{\gamma}^{(\alpha )}=\frac{\alpha ^{\left\vert \tilde{\tau}\right\vert
}}{\alpha ^{t}}\prod_{j=1}^{r}\frac{\Gamma \left( \left(
p\!+\!1\!-\!j\right) \!/\!\alpha \right) }{\Gamma \left( 1\!+\!\left(
r\!-\!j\right) \!/\!\alpha \right) }\frac{\prod_{s\in \tau }\left(
p\!+\!\alpha a^{\prime }\left( s\right) \!-\!l^{\prime }(s)\right) }{%
\prod_{s\in \tilde{\tau}}\left( r\!+\!\left( a^{\prime }\left( s\right)
\!+\!1\right) \alpha \!-\!l^{\prime }\left( s\right) \!-\!1\right) }.
\end{equation*}

Now consider the last ratio of the products in the above expression. For the
product term in the numerator that corresponds to square $s$ in the position 
$\left( i,j\right) $ in the diagram of $\tau ,$ there exists exactly the
same term in the denominator, which corresponds to square $s$ in the
position $(i,j+\left( p-r+1\right) /\alpha -1)$ in the diagram of $\tilde{%
\tau}$. Therefore, we can write%
\begin{equation*}
\tilde{\gamma}^{(\alpha )}=\frac{\alpha ^{\left\vert \tilde{\tau}\right\vert
}}{\alpha ^{t}}\prod_{j=1}^{r}\frac{\Gamma \left( \left(
p\!+\!1\!-\!j\right) \!/\!\alpha \right) }{\Gamma \left( 1\!+\!\left(
r\!-\!j\right) \!/\!\alpha \right) }\frac{1}{\prod_{s\in \hat{\tau}}\left(
r\!+\!\left( a^{\prime }\left( s\right) \!+\!1\right) \alpha \!-\!l^{\prime
}\left( s\right) \!-\!1\right) },
\end{equation*}%
where $\hat{\tau}$ is the partition that consists of $r$ identical parts $%
\left( p-r+1\right) /\alpha -1.$

Finally, note that%
\begin{eqnarray*}
&&\prod_{s\in \hat{\tau}}\left( r\!+\!\left( a^{\prime }\left( s\right)
\!+\!1\right) \alpha \!-\!l^{\prime }\left( s\right) \!-\!1\right) =\alpha
^{\left\vert \tilde{\tau}\right\vert -t}\prod_{s\in \hat{\tau}}\left( \left(
r\!-\!l^{\prime }(s)\!-1\right) /\alpha \!+\!a^{\prime }\left( s\right)
\!+\!1\right) \\
&=&\alpha ^{\left\vert \tilde{\tau}\right\vert -t}\prod_{j=1}^{r}\frac{%
\Gamma \left( \left( r\!-\!j\right) /\alpha \!+\!\left( p\!-\!r\!+\!1\right)
/\alpha \right) }{\Gamma \left( \left( r\!-\!j\right) /\alpha \!+\!1\right) }%
=\alpha ^{\left\vert \tilde{\tau}\right\vert -t}\prod_{j=1}^{r}\frac{\Gamma
\left( \left( p\!-\!j\!+\!1\right) /\alpha \right) }{\Gamma \left( \left(
r\!-\!j\right) /\alpha \!+\!1\right) }.
\end{eqnarray*}%
Therefore, $\tilde{\gamma}^{(\alpha )}\!=\!1$ and the statement of the lemma
follows from (\ref{F00}) and~(\ref{semi-final}).$\square \medskip $

\textsc{Proof of Proposition 2}

Proposition 1 and Corollary 1 directly imply (\ref{LR contour c 1}) and the
following formula for $L\left( h;\mu \right) $

\begin{equation}
L\left( h;\mu \right) =k_{1}\frac{n^{np}S^{pr-r\left( 1+r\right) /2}}{\Gamma
\left( np\right) }\int_{0}^{\infty }y^{p\left( n-r\right)
+r(r+1)/2-1}e^{-ny}\det R\mathrm{d}y,  \label{proof of P3}
\end{equation}%
where $R$ is an $r\times r$ matrix with%
\begin{equation*}
R_{ij}=\frac{1}{2\pi \mathrm{i}}\oint_{\mathcal{K}}e^{\frac{y}{S}\frac{nh_{i}%
}{1+h_{i}}z}z^{j-1}\prod_{s=1}^{p}\left( z-\lambda _{s}\right) ^{-1}\mathrm{d%
}z.
\end{equation*}

Let us write $\det R$ as%
\begin{equation*}
\det R=\sum_{\rho }\limfunc{sgn}\rho \prod_{j=1}^{r}\frac{1}{2\pi \mathrm{i}}%
\oint_{\mathcal{K}}e^{\frac{y}{S}\frac{nh_{\rho (j)}}{1+h_{\rho (j)}}%
z}z^{j-1}\prod_{s=1}^{p}\left( z-\lambda _{s}\right) ^{-1}\mathrm{d}z,
\end{equation*}%
or equivalently as 
\begin{equation}
\det R=\sum_{\rho }\frac{\limfunc{sgn}\rho }{\left( 2\pi \mathrm{i}\right)
^{r}}\oint_{\mathcal{K}}...\oint_{\mathcal{K}}\prod_{j=1}^{r}\left\{ e^{%
\frac{y}{S}\frac{nh_{\rho (j)}}{1+h_{\rho (j)}}z_{j}}z_{j}^{j-1}%
\prod_{s=1}^{p}\left( z_{j}-\lambda _{s}\right) ^{-1}\right\} \mathrm{d}%
z_{r}...\mathrm{d}z_{1}\text{.}  \label{representation of detR}
\end{equation}%
Using this representation, we have%
\begin{eqnarray*}
&&\int_{0}^{\infty }y^{p\left( n-r\right) +r(r+1)/2-1}e^{-ny}\det R(y)%
\mathrm{d}y=\sum_{\rho }\frac{\limfunc{sgn}\rho }{\left( 2\pi \mathrm{i}%
\right) ^{r}}\times \\
&&\int_{0}^{\infty }\oint_{\mathcal{K}}...\oint_{\mathcal{K}}y^{p\left(
n-r\right) +r(r+1)/2-1}\exp \left\{ -\left( n-\sum_{j=1}^{r}\frac{nh_{\rho
(j)}}{1+h_{\rho (j)}}\frac{z_{j}}{S}\right) y\right\} \times \\
&&\prod_{j=1}^{r}\left\{ z_{j}^{j-1}\prod_{s=1}^{p}\left( z_{j}-\lambda
_{s}\right) ^{-1}\right\} \mathrm{d}z_{r}...\mathrm{d}z_{1}\mathrm{d}y.
\end{eqnarray*}%
Since the contour $\mathcal{K}$ is chosen so that for any $z\in \mathcal{K}$%
, $\func{Re}z<\left( \sum_{j=1}^{r}\frac{h_{j}}{1+h_{j}}\right) ^{-1}S,$ the
integrand in the above multiple integral is absolutely integrable on $\left[
0,\infty \right) \times \mathcal{K}\times ...\times \mathcal{K}$, and
Fubini's theorem justifies the interchange of the order of the integrals, so
that%
\begin{eqnarray*}
&&\int_{0}^{\infty }y^{p\left( n-r\right) +r(r+1)/2-1}e^{-ny}\det R(y)%
\mathrm{d}y=\frac{\Gamma \left( p\left( n-r\right) +r(r+1)/2\right) }{%
n^{p\left( n-r\right) +r(r+1)/2}}\times \\
&&\sum_{\rho }\frac{\limfunc{sgn}\rho }{\left( 2\pi \mathrm{i}\right) ^{r}}%
\oint_{\mathcal{K}}...\oint_{\mathcal{K}}\left( 1-\sum_{j=1}^{r}\frac{%
h_{\rho (j)}}{1+h_{\rho (j)}}\frac{z_{j}}{S}\right) ^{-p\left( n-r\right)
-r(r+1)/2}\times \\
&&\prod_{j=1}^{r}\left\{ z_{j}^{j-1}\prod_{s=1}^{p}\left( z_{j}-\lambda
_{s}\right) ^{-1}\right\} \mathrm{d}z_{r}...\mathrm{d}z_{1}\text{.}
\end{eqnarray*}%
Combining this with (\ref{proof of P3}), we get (\ref{LR contour c 2}). $%
\square $\medskip

\textsc{Proof of Lemma 2}

The lemma can be proven using arguments very similar to those in the proof
of Lemmas 4 and 6 in Onatski, Moreira and Hallin (2012) (OMH in what
follows), and we omit the proof to save space.$\square $

\textsc{Proof of Lemma 3.}

To save space, we will only establish (\ref{Watson}), relegating a
conceptually similar but more technical proof of (\ref{Watson1}) to the
Supplementary Appendix. Lemma 5 in OMH implies that 
\begin{equation}
\oint_{\mathcal{K}_{i}}e^{-nf_{i}\left( z\right) }g\left( z\right) \mathrm{d}%
z=e^{-nf_{i0}}\left[ \frac{g\left( z_{i0}\right) \pi ^{1/2}}{%
f_{i2}^{1/2}n^{1/2}}+\frac{O_{p}\left( 1\right) }{h_{i}n^{3/2}}\right] ,
\label{OMH_lemma5}
\end{equation}%
where $g\left( z\right) =\exp \left\{ -\frac{1}{2}\Delta _{p}\left( z\right)
\right\} $ and $O_{p}\left( 1\right) $ is uniform in $h_{i}\in \left( 0,\bar{%
h}\right] $. A careful inspection of OMH's proof of their Lemma 5 reveals
that a version of (\ref{OMH_lemma5}) remains valid for general functions $%
g\left( z\right) $ that are analytic in the open ball $B\left(
z_{i0},r_{i}\right) $ with center at $z_{i0}$ and radius $r_{i}=\min \left\{
z_{i0}-\max \left\{ \bar{b}_{p},\lambda _{1}\right\} ,\frac{1+h_{i}}{h_{i}}%
S-z_{i0}\right\} $ with probability approaching 1 as $n,p\rightarrow \infty $%
. Precisely, for such general $g\left( z\right) $ we have 
\begin{equation}
\oint_{\mathcal{K}_{i}}e^{-nf_{i}\left( z\right) }g\left( z\right) \mathrm{d}%
z=e^{-nf_{i0}}\frac{g\left( z_{i0}\right) \pi ^{1/2}}{f_{i2}^{1/2}n^{1/2}}%
+\Psi _{1}+\Psi _{2}+\Psi _{3}  \label{OMH_relaxed}
\end{equation}%
with%
\begin{eqnarray}
\left\vert \Psi _{1}\right\vert
&<&C_{1}e^{-nf_{i0}}h_{i}^{-1}n^{-3/2}\sup_{z\in \bar{B}}\left\vert g\left(
z\right) \right\vert ,  \label{psi1} \\
\left\vert \Psi _{2}\right\vert
&<&C_{1}e^{-nf_{i0}}e^{-nC_{2}}h_{i}^{-1}\sup_{z\in \mathcal{K}_{i1}\cup 
\mathcal{\bar{K}}_{i1}}\left\vert g\left( z\right) \right\vert ,\text{ and}
\label{psi2} \\
\left\vert \Psi _{3}\right\vert &<&C_{1}\left\vert \oint_{\mathcal{K}%
_{i2}\cup \mathcal{\bar{K}}_{i2}}e^{-nf_{i}\left( z\right) }g\left( z\right) 
\mathrm{d}z\right\vert ,  \label{psi3}
\end{eqnarray}%
where $C_{1}$ and $C_{2}$ are some positive constants, and $\bar{B}$ is a
closed ball with center at $z_{i0}$ and radius $r_{i}/2$.

Now, let $g\left( z\right) =g_{j}(z)=z^{j-1}\exp \left\{ -\Delta _{p}\left(
z\right) \right\} $. Lemma A2 in OMH implies that $\sup_{z\in \bar{B}\cup 
\mathcal{K}_{i1}\cup \mathcal{\bar{K}}_{i1}}\left\vert g\left( z\right)
\right\vert =h_{i}^{1-j}O_{p}\left( 1\right) $ uniformly in $h_{i}\in \left(
0,\bar{h}\right] $. Therefore, by (\ref{psi1}) and (\ref{psi2}),%
\begin{equation}
\Psi _{1}+\Psi _{2}=e^{-nf_{i0}}h_{i}^{-j}n^{-3/2}O_{p}\left( 1\right) \text{%
.}  \label{psi1psi2}
\end{equation}

Turning to the analysis of $\Psi _{3}$, note that by definition of $%
f_{i}\left( z\right) $ and $g\left( z\right) ,$%
\begin{equation}
e^{-nf_{i}\left( z\right) }g\left( z\right) =e^{n\frac{h_{i}}{1+h_{i}}%
z}z^{j-1}\prod_{j=1}^{p}\left( z-\lambda _{j}\right) ^{-1}.  \label{fact}
\end{equation}%
For $z\in \mathcal{K}_{i2}\cup \mathcal{\bar{K}}_{i2},$ we have $\left\vert
\left( z-\lambda _{j}\right) ^{-1}\right\vert <\left( 3z_{i0}\right) ^{-1},$
and $\left\vert z\left( z-\lambda _{j}\right) ^{-1}\right\vert <2,$ for any $%
j=1,...,p$. Therefore, using (\ref{fact}), we get%
\begin{eqnarray*}
\left\vert \oint_{\mathcal{K}_{i2}\cup \mathcal{\bar{K}}_{i2}}e^{-nf_{i}%
\left( z\right) }g\left( z\right) \mathrm{d}z\right\vert &<&2^{j-1}\left(
3z_{i0}\right) ^{-p+j-1}\oint_{\mathcal{K}_{i2}\cup \mathcal{\bar{K}}%
_{i2}}\left\vert e^{n\frac{h_{i}}{1+h_{i}}z}\mathrm{d}z\right\vert \\
&=&2^{j}\left( 3z_{i0}\right) ^{-p+j-1}\left( n\frac{h_{i}}{1+h_{i}}\right)
^{-1}e^{n\frac{h_{i}}{1+h_{i}}z_{i0}} \\
&=&2^{j}\left( 3z_{i0}\right) ^{j-1}\left( n\frac{h_{i}}{1+h_{i}}\right)
^{-1}e^{-n\left( c_{p}\ln \left( 3z_{i0}\right) -\frac{h_{i}}{1+h_{i}}%
z_{i0}\right) } \\
&=&2^{j}\left( 3z_{i0}\right) ^{j-1}\left( n\frac{h_{i}}{1+h_{i}}\right)
^{-1}3^{-p}e^{-n\left( c_{p}\ln \left( z_{i0}\right) -h_{i}-c_{p}\right) }.
\end{eqnarray*}%
On the other hand, for any $h_{i}\in \left[ 0,\bar{h}\right] ,$ $h_{i}<\sqrt{%
c_{p}}$ for sufficiently large $n$ and $p,$ and $c_{p}\ln \left(
z_{i0}\right) -h_{i}-c_{p}>f_{i0}.$ Indeed, using the definition of $z_{i0}$
and the fact, established in OMH's Lemma 11, that $f_{i0}=-c_{p}-\left(
1-c_{p}\right) \ln \left( 1+h_{i}\right) +c_{p}\ln \frac{c_{p}}{h_{i}},$ we
have%
\begin{equation*}
c_{p}\ln \left( z_{i0}\right) -h_{i}-c_{p}-f_{i0}=\ln \left( 1+h_{i}\right)
+c_{p}\ln \left( c_{p}+h_{i}\right) -h_{i}-c_{p}\ln c_{p}.
\end{equation*}%
The right hand side of this equality equals $0$ at $h_{i}=0$ and has a
non-negative derivative with respect to $h_{i}$ for all $0\leq h_{i}\leq 
\sqrt{c_{p}}.$ Therefore, 
\begin{equation*}
\left\vert \oint_{\mathcal{K}_{i2}\cup \mathcal{\bar{K}}_{i2}}e^{-nf_{i}%
\left( z\right) }g\left( z\right) \mathrm{d}z\right\vert <2^{j}\left(
3z_{i0}\right) ^{j-1}\left( n\frac{h_{i}}{1+h_{i}}\right)
^{-1}3^{-p}e^{-nf_{i0}},
\end{equation*}%
and thus, $\Psi _{3}=e^{-nf_{i0}}h_{i}^{-j}n^{-3/2}O_{p}\left( 1\right) $,
uniformly in $h_{i}\in \left( 0,\bar{h}\right] $. Combining this with (\ref%
{OMH_relaxed}) and (\ref{psi1psi2}), we obtain (\ref{Watson}).$\square $%
\medskip

\textsc{Proof of Theorem 1}

Proposition 2 and Lemma 3 imply that 
\begin{equation*}
L\left( h;\lambda \right) =\frac{k_{1}\exp \left\{
-n\sum_{i=1}^{r}f_{i0}\right\} }{\left( 2\mathrm{i}\right) ^{r}\left( \pi
n\right) ^{r/2}}\det \left( \frac{z_{i0}^{j-1}\exp \left\{ -\Delta
_{p}\left( z_{i0}\right) \right\} }{f_{i2}^{1/2}}+\frac{O_{p}\left( 1\right) 
}{h_{i}^{j}n}\right) _{1\leq i,j\leq r}
\end{equation*}%
As is shown in OMH (see their Lemma 11 and (A8))\footnote{%
Note that the expressions given in OMH are half times the expressions given
below because the equivalent of $f_{i}$ in the real-valued data case
considered by OMH is $f_{i}/2$.}, for $h_{i}\leq \bar{h}$, 
\begin{eqnarray}
f_{i0} &=&-c_{p}-\left( 1-c_{p}\right) \ln \left( 1+h_{i}\right) +c_{p}\ln 
\frac{c_{p}}{h_{i}},\text{ and}  \label{fi0} \\
f_{i2} &=&-\frac{h_{i}^{2}}{2\left( 1+h_{i}\right) ^{2}\left(
c_{p}-h_{i}^{2}\right) }.  \label{fi2}
\end{eqnarray}%
Moreover, by OMH's Lemma A2, $\exp \left\{ -\Delta _{p}\left( z_{i0}\right)
\right\} =O_{p}\left( 1\right) $ uniformly in $h\in \left( 0,\bar{h}\right]
^{r}$. Using these facts and the definition of $k_{1}$ given in Proposition
2, we get after some algebra

\begin{eqnarray*}
L\left( h;\lambda \right) &=&n^{r^{2}/2}\prod_{t=1}^{r}\left[ \left(
p-t\right) !\left( \frac{c_{p}-h_{t}^{2}}{2\pi }\right) ^{1/2}\right]
e^{rp}p^{-rp}\times \\
&&\exp \left\{ -\sum_{i=1}^{r}\Delta _{p}\left( z_{i0}\right) \right\}
\prod_{i>j}^{r}\left( c_{p}-h_{i}h_{j}\right) \left( 1+O_{p}\left(
n^{-1}\right) \right)
\end{eqnarray*}%
Applying Stirling's formula 
\begin{equation*}
\left( p-t\right) !=e^{-p}p^{p-t+1}\left( \frac{2\pi }{p}\right)
^{1/2}\left( 1+O\left( p^{-1}\right) \right)
\end{equation*}%
we get%
\begin{equation*}
L\left( h;\lambda \right) =\exp \left\{ -\sum_{i=1}^{r}\Delta _{p}\left(
z_{i0}\right) \right\} \prod_{t=1}^{r}\left( 1-\frac{h_{t}^{2}}{c_{p}}%
\right) ^{1/2}\prod_{i>j}^{r}\left( 1-\frac{h_{i}h_{j}}{c_{p}}\right) \left(
1+O_{p}\left( n^{-1}\right) \right) ,
\end{equation*}%
which implies (\ref{equivalence 1}).

Turning to the proof of (\ref{equivalence 2}), Proposition 2 and Lemma 3
imply that 
\begin{eqnarray}
&&L\left( h;\mu \right) =\left( -1\right) ^{r\left( r-1\right)
/2}n^{-pr+r\left( r+1\right) /2}\prod_{i>j}^{r}\left( h_{i}-h_{j}\right)
^{-1}\times  \label{Lhm} \\
&&\prod_{t=1}^{r}\left[ h_{t}^{r-p}\left( 1+h_{t}\right) ^{p-n-1}\left(
p-t\right) !\right] \left( nS\right) ^{pr-r(r+1)/2}\frac{\Gamma \left(
p\left( n-r\right) +r(r+1)/2\right) }{\Gamma \left( np\right) }\times  \notag
\\
&&\sum_{\rho }\frac{\limfunc{sgn}\rho }{\left( 2\pi \mathrm{i}\right) ^{r}}%
q_{\rho }\left( \mathbf{z}_{0}\right) \prod_{j=1}^{r}e^{-nf_{\rho \left(
j\right) 0}}\frac{g_{j}\left( z_{\rho \left( j\right) 0}\right) \pi ^{1/2}}{%
f_{\rho \left( j\right) 2}^{1/2}n^{1/2}}\left( 1+O_{p}\left( n^{-1}\right)
\right)  \notag
\end{eqnarray}%
Using the definition of $q_{\rho }\left( \mathbf{z}_{0}\right) $ and of $%
z_{i0},$ we get%
\begin{eqnarray}
q_{\rho }\left( \mathbf{z}_{0}\right) &=&\left( 1-\sum_{i=1}^{r}\frac{h_{i}}{%
1+h_{i}}\frac{z_{i0}}{S}\right) ^{-p\left( n-r\right) -r(r+1)/2}\exp \left\{
-\sum_{i=1}^{r}\frac{nh_{i}z_{i0}}{1+h_{i}}\right\}  \notag \\
&=&\left( 1-\sum_{i=1}^{r}\frac{h_{i}+c_{p}}{S}\right) ^{-p\left( n-r\right)
-r(r+1)/2}\exp \left\{ -pr-\sum_{i=1}^{r}nh_{i}\right\} .  \label{qpozo}
\end{eqnarray}%
Further, using the definition of $g_{j}\left( z_{\rho \left( j\right)
0}\right) $, the fact that $\sum_{\rho }\limfunc{sgn}\rho z_{\rho \left(
j\right) 0}^{j-1}$ equals the Vandermonde determinant $\prod_{i>j}^{r}\left(
z_{i0}-z_{j0}\right) ,$ we get%
\begin{eqnarray}
\sum_{\rho }\limfunc{sgn}\rho \prod_{j=1}^{r}\frac{e^{-nf_{\rho \left(
j\right) 0}}g_{j}\left( z_{\rho \left( j\right) 0}\right) \pi ^{1/2}}{%
f_{\rho \left( j\right) 2}^{1/2}n^{1/2}}\!\! &=&\!\!\frac{\pi ^{r/2}}{n^{r/2}%
}\exp \left\{ \!-\!\sum\limits_{i=1}^{r}\left( nf_{i0}+\Delta _{p}\left(
z_{i0}\right) \right) \right\} \!\!\times  \label{sgnroetc} \\
&&\prod_{i=1}^{r}f_{i2}^{-1/2}\prod_{i>j}^{r}\left( z_{i0}-z_{j0}\right) . 
\notag
\end{eqnarray}%
Substituting (\ref{qpozo}) and (\ref{sgnroetc}) into (\ref{Lhm}), and using (%
\ref{fi0}) and (\ref{fi2}) together with the fact that the branch of the
square root in $f_{i2}^{-1/2}$ is chosen so that $\sqrt{-1}=-\mathrm{i},$ we
get after some algebra 
\begin{eqnarray*}
L\left( h;\mu \right) &=&\prod_{t=1}^{r}\left[ \left( p-t\right) !\left(
c_{p}-h_{t}^{2}\right) ^{1/2}\right] S^{pr-r(r+1)/2}\times \\
&&\frac{\Gamma \left( p\left( n-r\right) +r(r+1)/2\right) }{\Gamma \left(
np\right) }\left( 1-\sum_{j=1}^{r}\frac{h_{j}+c_{p}}{S}\right) ^{-p\left(
n-r\right) -r(r+1)/2}\times \\
&&\exp \left\{ -\sum_{j=1}^{r}nh_{j}\right\} \frac{1}{\left( 2\pi n\right)
^{r/2}}c_{p}^{-pr}\prod_{i>j}^{r}\left( c_{p}-h_{i}h_{j}\right) \left(
1+O_{p}\left( n^{-1}\right) \right) .
\end{eqnarray*}%
Now, using the fact that $S-p=O_{p}\left( 1\right) ,$ we get $\ln \left( 
\frac{S}{p}\right) =\frac{S-p}{p}+O_{p}\left( p^{-2}\right) $ and 
\begin{eqnarray*}
\ln \left( 1-\sum_{j=1}^{r}\frac{h_{j}+c_{p}}{S}\right) &=&-\frac{%
\sum_{j=1}^{r}\left( h_{j}+c_{p}\right) }{p}-\frac{1}{2}\frac{\left(
\sum_{j=1}^{r}\left( h_{j}+c_{p}\right) \right) ^{2}}{p^{2}}+ \\
&&\frac{\sum_{j=1}^{r}\left( h_{j}+c_{p}\right) }{p^{2}}\left( S-p\right)
+O_{p}\left( p^{-3}\right) .
\end{eqnarray*}%
Further, the Stirling approximations give%
\begin{eqnarray*}
\left( p-t\right) ! &=&e^{-p}p^{p-t+1}\left( \frac{2\pi }{p}\right)
^{1/2}\left( 1+O\left( p^{-1}\right) \right) \text{ and} \\
\frac{\Gamma \left( p\left( n-r\right) +r(r+1)/2\right) }{\Gamma \left(
np\right) } &=&\left( pn\right) ^{-pr+r\left( r+1\right) /2}e^{\frac{1}{2}%
c_{p}r^{2}}\left( 1+O\left( n^{-1}\right) \right) .
\end{eqnarray*}%
So finally, after some cancellations,%
\begin{eqnarray*}
L\left( h;\mu \right) &=&e^{\frac{1}{2c_{p}}\left(
\sum_{j=1}^{r}h_{j}\right) ^{2}}e^{-\frac{\sum_{j=1}^{r}h_{j}}{c_{p}}\left(
S-p\right) }e^{-\sum\limits_{i=1}^{r}\Delta _{p}\left( z_{i0}\right) } \\
&&\prod_{t=1}^{r}\left( 1-\frac{h_{t}^{2}}{c_{p}}\right)
^{1/2}\prod_{i>j}^{r}\left( 1-\frac{h_{i}h_{j}}{c_{p}}\right) \left(
1+O_{p}\left( n^{-1}\right) \right) ,
\end{eqnarray*}%
which implies (\ref{equivalence 2}).

To establish the rest of the statements of Theorem 1 we will need the
following lemma.\medskip

\textbf{Lemma A3. }\textit{Suppose that our null hypothesis holds. Denote }$%
\sum_{j=1}^{p}\lambda _{j}^{2}$\textit{\ as }$T.$\textit{\ Then, for any
fixed }$r$\textit{\ and }$\bar{h}<\sqrt{c},$\textit{\ and any }$\left(
h_{1},...,h_{r}\right) \in \left( 0,\bar{h}\right] ^{r},$\textit{\ as }$%
n,p\rightarrow \infty $\textit{\ so that }$p/n\rightarrow c$\textit{, the
vector }$\left( S-p,T-\left( 1+c_{p}\right) p,\Delta _{p}\left(
z_{10}\right) ,...,\Delta _{p}\left( z_{r0}\right) \right) $\textit{\
converges in distribution to a Gaussian vector }$\left( \eta ,\zeta ,\xi
_{1},...,\xi _{r}\right) $\textit{\ with }%
\begin{eqnarray*}
\mathrm{E}\left( \eta \right) &=&\mathrm{E}\left( \zeta \right) =\mathrm{E}%
\left( \xi _{i}\right) =0,\text{ } \\
\mathrm{Var}\left( \eta \right) &=&c,\mathrm{Var}\left( \zeta \right)
=2c\left( 2+5c+2c^{2}\right) ,\mathrm{Cov}\left( \eta ,\zeta \right)
=2c\left( 1+c\right) , \\
\mathrm{Cov}\left( \eta ,\xi _{i}\right) &=&-h_{i},\mathrm{Cov}\left( \zeta
,\xi _{i}\right) =-h_{i}\left( h_{i}+2+2c\right) ,\text{ and} \\
\mathrm{Cov}\left( \xi _{i},\xi _{k}\right) &=&-\ln \left(
1-h_{i}h_{k}/c\right)
\end{eqnarray*}

\textsc{Proof: }The proof of the lemma is similar to that of Lemma 12 in
OMH. The convergence to the Gaussian distribution follows from Theorem 1.1
of Bai and Silverstein (2004). The formulas for the means, variances and
covariances of $\eta $ and $\xi _{j}$ are obtained using Theorem 1.1 iii) of
Bai and Silverstein (2004) similarly to how the corresponding formulas in
Lemma 12 of OMH are obtained using Theorem 1.1 ii). Therefore, below we only
derive the formulae for the mean, variance, and covariances that involve $%
\zeta $. Variable $\zeta $ does not appear in Lemma 12 of OMH because the
lemma does not study the asymptotics of $T-\left( 1+c_{p}\right) p.$

The fact that $\mathrm{E}\zeta =0$ follows directly from Theorem 1.1 iii) of
Bai and Silverstein (2004). The same theorem implies that 
\begin{equation}
\limfunc{Cov}\left( \xi _{j},\zeta \right) =-\frac{1}{4\pi ^{2}}\oint \oint 
\frac{z_{2}^{2}\ln \left( \bar{z}_{j0}-z_{1}\right) }{\left( \underline{m}%
\left( z_{1}\right) -\underline{m}\left( z_{2}\right) \right) ^{2}}\frac{%
\mathrm{d}\underline{m}\left( z_{1}\right) }{\mathrm{d}z_{1}}\frac{\mathrm{d}%
\underline{m}\left( z_{2}\right) }{\mathrm{d}z_{2}}\mathrm{d}z_{1}\mathrm{d}%
z_{2},  \label{dzeta_ksi}
\end{equation}%
where $\bar{z}_{j0}=\lim z_{j0}$ as $n,p\rightarrow \infty ,$%
\begin{equation}
\limfunc{Cov}\left( \eta ,\zeta \right) =-\frac{1}{4\pi ^{2}}\oint \oint 
\frac{z_{2}^{2}z_{1}}{\left( \underline{m}\left( z_{1}\right) -\underline{m}%
\left( z_{2}\right) \right) ^{2}}\frac{\mathrm{d}\underline{m}\left(
z_{1}\right) }{\mathrm{d}z_{1}}\frac{\mathrm{d}\underline{m}\left(
z_{2}\right) }{\mathrm{d}z_{2}}\mathrm{d}z_{1}\mathrm{d}z_{2},
\label{eta_dzeta}
\end{equation}%
and 
\begin{equation}
\limfunc{Var}\left( \zeta \right) =-\frac{1}{4\pi ^{2}}\oint \oint \frac{%
z_{1}^{2}z_{2}^{2}}{\left( \underline{m}\left( z_{1}\right) -\underline{m}%
\left( z_{2}\right) \right) ^{2}}\frac{\mathrm{d}\underline{m}\left(
z_{1}\right) }{\mathrm{d}z_{1}}\frac{\mathrm{d}\underline{m}\left(
z_{2}\right) }{\mathrm{d}z_{2}}\mathrm{d}z_{1}\mathrm{d}z_{2},
\label{dzeta_variance}
\end{equation}%
where 
\begin{equation*}
\underline{m}\left( z\right) =-\left( 1-c\right) z^{-1}+cm(z)
\end{equation*}%
with $m\left( z\right) $ given by (3.6) of OMH, where $c_{p}$ is replaced by 
$c.$ That is,%
\begin{equation}
\underline{m}\left( z\right) =\frac{-z+c-1+\sqrt{\left( z-c-1\right) ^{2}-4c}%
}{2z},  \label{m lower bar 1}
\end{equation}%
where the branch of the square root is chosen so that the real and the
imaginary parts of $\sqrt{\left( z-c-1\right) ^{2}-4c}$ have the same signs
as the real and the imaginary parts of $z-c-1$, respectively. The contours
of integration in (\ref{dzeta_ksi})-(\ref{dzeta_variance}) are closed,
oriented counterclockwise, enclose zero and the support of the
Marchenko-Pastur distribution with parameter $c$, and do not enclose $\bar{z}%
_{j0}$.

The above expressions can be simplified. Use formula 1.16 of Bai and
Silverstein (2004), to get%
\begin{equation}
\limfunc{Cov}\left( \xi _{j},\zeta \right) =-\frac{1}{4\pi ^{2}}\oint \oint 
\frac{\ln \left( \bar{z}_{j0}-z\left( m_{1}\right) \right) \left(
z(m_{2})\right) ^{2}}{\left( m_{1}-m_{2}\right) ^{2}}\mathrm{d}m_{1}\mathrm{d%
}m_{2},  \label{simplified cov}
\end{equation}%
\begin{equation}
\limfunc{Cov}\left( \eta ,\zeta \right) =-\frac{1}{4\pi ^{2}}\oint \oint 
\frac{z\left( m_{1}\right) \left( z(m_{2})\right) ^{2}}{\left(
m_{1}-m_{2}\right) ^{2}}\mathrm{d}m_{1}\mathrm{d}m_{2},\text{ and}
\label{simplified_cov_eta_dzeta}
\end{equation}%
\begin{equation}
\limfunc{Var}\left( \zeta \right) =-\frac{1}{4\pi ^{2}}\oint \oint \frac{%
\left( z\left( m_{1}\right) \right) ^{2}\left( z(m_{2})\right) ^{2}}{\left(
m_{1}-m_{2}\right) ^{2}}\mathrm{d}m_{1}\mathrm{d}m_{2},
\label{Simplified_var_dzeta}
\end{equation}%
where 
\begin{equation}
z\left( m\right) =-\frac{1}{m}+\frac{c}{1+m}  \label{inverse Stiltjes}
\end{equation}%
and the contours of integration over $m_{1}$ and $m_{2}$ in (\ref{simplified
cov}-\ref{Simplified_var_dzeta}) are obtained from the contours of
integration over $z_{1}$ and $z_{2}$ in (\ref{dzeta_ksi}-\ref{dzeta_variance}%
) by transformation $\underline{m}\left( z\right) .$ Recall that by
assumption the contours over $z_{1}$ and $z_{2}$ intersect the real line to
the left of zero and in between the upper boundary of the support of the
Marchenko-Pastur distribution, $\left( 1+\sqrt{c}\right) ^{2}$, and $\bar{z}%
_{j0}.$ Therefore, as can be shown using the definition (\ref{m lower bar 1}%
) of $\underline{m}\left( z\right) $, the $m_{1}$-contour and $m_{2}$%
-contour are clockwise oriented and intersect the real line in between $%
-\left( 1+\sqrt{c}\right) ^{-1}$ and $\underline{m}\left( \bar{z}%
_{j0}\right) =-h_{j}\left( h_{j}+c\right) ^{-1}$ and to the right of zero.
In particular, both contours enclose $0$ and $-h_{j}\left( h_{j}+c\right)
^{-1}$, but not $-1$ and $-\left( 1+h_{j}\right) ^{-1}$.

Assuming without loss of generality that $m_{1}$-contour lies inside the $%
m_{2}$-contour, from (A64) in the Supplementary appendix of OMH, we have%
\begin{equation}
\oint \frac{\ln \left( \bar{z}_{j0}-z\left( m_{1}\right) \right) }{\left(
m_{1}-m_{2}\right) ^{2}}\mathrm{d}m_{1}=2\pi \mathrm{i}\left( -\frac{1}{m_{2}%
}+\frac{1}{m_{2}+h_{j}\left( h_{j}+c\right) ^{-1}}\right) .
\label{inner contour}
\end{equation}%
Denoting $-h_{j}\left( h_{j}+c\right) ^{-1}$ as $x_{j},$ we get from (\ref%
{inner contour}) and (\ref{simplified cov})%
\begin{eqnarray*}
\limfunc{Cov}\left( \xi _{j},\zeta \right) &=&\frac{2\pi \mathrm{i}}{4\pi
^{2}}\oint \left( z(m_{2})\right) ^{2}\left( \frac{1}{m_{2}}-\frac{1}{%
m_{2}-x_{j}}\right) \mathrm{d}m_{2} \\
&=&\frac{2\pi \mathrm{i}}{4\pi ^{2}}\oint \left( -\frac{1}{m_{2}}+\frac{c}{%
1+m_{2}}\right) ^{2}\left( \frac{1}{m_{2}}-\frac{1}{m_{2}-x_{j}}\right) 
\mathrm{d}m_{2} \\
&=&-h_{j}\left( h_{j}+2+2c\right) ,
\end{eqnarray*}%
where the last equality follows from Cauchy's residue theorem and the fact
that the contour is oriented clock-wise.

For $\limfunc{Cov}\left( \eta ,\zeta \right) ,$ we have%
\begin{equation*}
\oint \frac{z\left( m_{1}\right) }{\left( m_{1}-m_{2}\right) ^{2}}\mathrm{d}%
m_{1}=\oint \frac{\left( -\frac{1}{m_{1}}+\frac{c}{1+m_{1}}\right) }{\left(
m_{1}-m_{2}\right) ^{2}}\mathrm{d}m_{1}=\frac{2\pi \mathrm{i}}{m_{2}^{2}}
\end{equation*}%
so that%
\begin{eqnarray*}
\limfunc{Cov}\left( \eta ,\zeta \right) &=&-\frac{2\pi \mathrm{i}}{4\pi ^{2}}%
\oint \oint \frac{\left( z(m_{2})\right) ^{2}}{m_{2}^{2}}\mathrm{d}m_{2} \\
&=&-\frac{2\pi \mathrm{i}}{4\pi ^{2}}\oint \left( -\frac{1}{m_{2}}+\frac{c}{%
1+m_{2}}\right) ^{2}\frac{1}{m_{2}^{2}}\mathrm{d}m_{2} \\
&=&2c\left( 1+c\right)
\end{eqnarray*}%
by Cauchy's theorem.

For $\limfunc{Var}\left( \zeta \right) ,$ we have%
\begin{eqnarray*}
\oint \frac{z\left( m_{1}\right) ^{2}}{\left( m_{1}-m_{2}\right) ^{2}}%
\mathrm{d}m_{1} &=&\oint \frac{\left( -\frac{1}{m_{1}}+\frac{c}{1+m_{1}}%
\right) ^{2}}{\left( m_{1}-m_{2}\right) ^{2}}\mathrm{d}m_{1} \\
&=&\frac{4\pi \mathrm{i}}{m_{2}^{2}}\left( c-\frac{1}{m_{2}}\right)
\end{eqnarray*}%
so that%
\begin{eqnarray*}
\limfunc{Var}\left( \zeta \right) &=&-\frac{4\pi \mathrm{i}}{4\pi ^{2}}\oint
\oint \frac{\left( z(m_{2})\right) ^{2}}{m_{2}^{2}}\left( c-\frac{1}{m_{2}}%
\right) \mathrm{d}m_{2} \\
&=&-\frac{4\pi \mathrm{i}}{4\pi ^{2}}\oint \oint \frac{\left( -\frac{1}{m_{2}%
}+\frac{c}{1+m_{2}}\right) ^{2}\left( c-\frac{1}{m_{2}}\right) }{m_{2}^{2}}%
\mathrm{d}m_{2} \\
&=&2c\left( 2+5c+2c^{2}\right)
\end{eqnarray*}%
by Cauchy's theorem.$\square $

Lemma A3 and formulae (\ref{equivalence 1}) and (\ref{equivalence 2}) imply
the convergence of finite dimensional-distributions of the random fields $%
\ln L\left( h;\lambda \right) $ and $\ln L\left( h;\mu \right) $ to the
Gaussian distributions with means and covariance matrices characterized by (%
\ref{mean}-\ref{covariance mu}).

To complete the proof of Theorem 1, we need to establish the tightness of $%
\ln L\left( h;\lambda \right) $ and $\ln L\left( h;\mu \right) $, viewed as
random elements of the space $C\left[ 0,\overline{h}\right] ^{r},$ as $%
n,p\rightarrow \infty $ so that $p/n\rightarrow c.$ Formulae (\ref%
{equivalence 1}-\ref{equivalence 2}) and the facts that $S-p=O_{p}\left(
1\right) ,$ and that $\Delta _{p}\left( z_{i0}\right) =O_{p}\left( 1\right) $
for $i=1,...,r,$ where $O_{p}\left( 1\right) $ are uniform in $h\in \left( 0,%
\bar{h}\right] ^{r},$ imply that for an arbitrarily small positive $%
\varepsilon ,$ there must exist $B>0$ such that $\Pr \left( \sup_{h\in
\left( 0,\bar{h}\right] ^{r}}\left\vert \ln L\left( h;\lambda \right)
\right\vert >B\right) <\varepsilon $ and $\Pr \left( \sup_{h\in \left( 0,%
\bar{h}\right] ^{r}}\left\vert \ln L\left( h;\mu \right) \right\vert
>B\right) <\varepsilon $ for sufficiently large $n$ and $p.$ Since, as
implied by Proposition 1, $\ln L\left( h;\lambda \right) $ and $\ln L\left(
h;\mu \right) $ are continuous functions on $h\in \left[ 0,\overline{h}%
\right] ^{r},$ $\sup_{h\in \left( 0,\bar{h}\right] ^{r}}\left\vert \ln
L\left( h;\lambda \right) \right\vert =\sup_{h\in \left[ 0,\bar{h}\right]
^{r}}\left\vert \ln L\left( h;\lambda \right) \right\vert ,$ and $\sup_{h\in
\left( 0,\bar{h}\right] ^{r}}\left\vert \ln L\left( h;\mu \right)
\right\vert =\sup_{h\in \left[ 0,\bar{h}\right] ^{r}}\left\vert \ln L\left(
h;\mu \right) \right\vert ,$ so that the tightness of $\ln L\left( h;\lambda
\right) $ and $\ln L\left( h;\mu \right) $ follows.$\square $

\textsc{Proof of theorem 2}

To save space, we only derive the asymptotic power envelope for the
relatively more difficult case of real-valued data and $\mu $-based tests.
According to the Neyman-Pearson lemma, the most powerful test of the null $%
h=0$ against a point alternative $h=\left( h_{1},...,h_{r}\right) $ is the
test which rejects the null when $\ln L\left( h;\mu \right) $ is larger than
a critical value $C.$ It follows from Theorem 1 that, for such a test to
have asymptotic size $\alpha $, $C$ must be 
\begin{equation}
C=\sqrt{W\left( h\right) }\Phi ^{-1}\left( 1-\alpha \right) +m\left(
h\right) ,  \label{ctitical value}
\end{equation}%
where 
\begin{eqnarray*}
m\left( h\right) &=&\frac{1}{2}\sum_{i,j=1}^{r}\left( \ln \left( 1-\frac{%
h_{i}h_{j}}{c}\right) +\frac{h_{i}h_{j}}{c}\right) \text{ and} \\
W\left( h\right) &=&-\sum_{i,j=1}^{r}\left( \ln \left( 1-\frac{h_{i}h_{j}}{c}%
\right) +\frac{h_{i}h_{j}}{c}\right) .
\end{eqnarray*}%
Now, according to Le Cam's third lemma and Theorem 1, under $h=\left(
h_{1},...,h_{r}\right) ,$ $\ln L\left( h;\mu \right) \overset{d}{\rightarrow 
}N\left( m\left( h\right) +W\left( h\right) ,W\left( h\right) \right) .$
Therefore, the asymptotic power $\beta _{\mu }\left( h\right) $ is (\ref%
{local power mu}).$\square $


\begin{thebibliography}{99}
\bibitem{A} Andreief, C. (1883). \textquotedblleft Note sur une relation les
int\'{e}grales d\'{e}finies des produits des fonctions\textquotedblright , M%
\'{e}m.de la Soc. Sci. Bordeaux 2.

\bibitem{BaiS/} Bai, Z.D. and J.W. Silverstein (2004) \textquotedblleft CLT
for Linear Spectral Statistics of Large-Dimensional Sample Covariance
Matrices\textquotedblright , \textit{Annals of Probability} 32, 553-605.

\bibitem{BBP/} Baik, J., Ben Arous, G. and S. P\'{e}ch\'{e}. (2005)
\textquotedblleft Phase transition of the largest eigenvalue for non-null
complex sample covariance matrices\textquotedblright\ \textit{Annals of
Probability} 33, 1643--1697.

\bibitem{BDMN} Bianchi, P., M. Debbah, M. Ma\"{\i}da and J. Najim (2010)
\textquotedblleft Performance of Statistical Tests for Single Source
Detection using Random Matrix Theory\textquotedblright , manuscript.

\bibitem{CS} Collins, B. and P. \'{S}niady (2007), \textquotedblleft New
scaling of Itzykson-Zuber integrals\textquotedblright , \textit{Annales de
l'IHP Prob. Stats.} 43 (2), 139--146.

\bibitem{DE} Dumitriu, I., Edelman, A., and G. Shuman (2007)
\textquotedblleft MOPS: Multivariate orthogonal polynomials
(symbolically)\textquotedblright , \textit{Journal of Symbolic Computation}
42, 587--620

\bibitem{F} Forrester, P.J. (2011), \textquotedblleft Probability densities
and distributions for spiked Wishart $\beta $-ensembles\textquotedblright ,
ArXiv:1101.2261.v1

\bibitem{G} Goodman, N. R. (1963) \textquotedblleft Statistical analysis
based on a certain multivariate complex Gaussian distribution, (An
introduction).\textquotedblright\ \textit{Annals of Mathematical Statistics}
34, 152-177.

\bibitem{GR} Gross, K.I., and Richards, D.S.P. (1987), \textquotedblleft
Special Functions of Matrix Argument.I: Algebraic Induction, Zonal
Polynomials, and Hypergeometric Functions\textquotedblright , \textit{%
Transactions of the American Mathematical Society} 301 (2), 781-811

\bibitem{GM} Guionnet, A. and M. Ma\"{\i}da (2005), \textquotedblleft A
Fourier view on the R-transform and related asymptotics of spherical
integrals\textquotedblright , \textit{Journal of Functional Analysis} 222
(2), 435--490.

\bibitem{H-C/} Harish-Chandra (1957) \textquotedblleft Differential
Operators on Semi-simple Lie Algebra\textquotedblright , \textit{American
journal of Mathematics} 79, 87-120.

\bibitem{IZ/} Itzykson, C., and Zuber, J.B. (1980) \textquotedblleft The
Planar Approximation. II\textquotedblright , \textit{Journal of Mathematical
Physics} 21, 411-421.

\bibitem{Ja/} James, A. T. (1964) \textquotedblleft Distributions of matrix
variates and latent roots derived from normal samples\textquotedblright , \ 
\textit{Annals of Mathematical Statistics} 35, 475-501.

\bibitem{J/} Johnstone, I.M. (2001) \textquotedblleft On the distribution of
the largest eigenvalue in principal components analysis.\textquotedblright\ 
\textit{Annals of Statistics} 29, 295--327.

\bibitem{J2007} Johnstone, I.M. (2007). \textquotedblleft High dimensional
statistical inference and random matrices\textquotedblright , Proceedings of
the International Congress of Mathematicians, Madrid, Spain, 2006. European
Mathematical Society, 307-333.

\bibitem{JT} Johnstone, I.M., and D.M. Titterington (2009).
\textquotedblleft Statistical challenges of high-dimensional
data\textquotedblright , \textit{Philosophical Transactions of Royal Society
A} 367, 4237-4253

\bibitem{KE} Koev, P. and A. Edelman (2006). \textquotedblleft The Efficient
Evaluation of the Hypergeometric Function of a matrix
Argument\textquotedblright , \textit{Mathematics of Computation} 75 (254),
833-846

\bibitem{KN} Krichman, S., and Nadler, B. (2009) \textquotedblleft
Non-Parametric Detection of the Number of Signals: Hypothesis Testing and
Random Matrix Theory\textquotedblright , \textit{IEEE Transactions on Signal
Processing }57, 3930-3941.

\bibitem{Mac/} Macdonald, I. G. (1995) Symmetric functions and Hall
polynomials. Second edition. Oxford Mathematical Monographs. Oxford Science
Publications. The Clarendon Press, Oxford University Press, New York.

\bibitem{Mo/} Mo, M.Y. (2011) \textquotedblleft The rank 1 real Wishart
spiked model\textquotedblright , arXiv:1101.5144v1

\bibitem{MP} Marchenko, V.A., and L.A. Pastur (1967) \textquotedblleft
Distribution of eigenvalues for some sets of random
matrices\textquotedblright , Math. USSR-Sbornik, vol. 1, no. 4, 457-483

\bibitem{Mar} Marinari, E., Parisi, G., and Ritort, F. (1994),
\textquotedblleft Replica field theory for deterministic models. II. A
non-random spin glass with glassy behavior\textquotedblright , \textit{J.
Phys. A} 27 (23), 7647--7668.

\bibitem{Mul} Muller, R.R., Guo, D., and Moustakas, A.L (2008)
\textquotedblleft Vector precoding for wireless MIMO systems and its replica
analysis\textquotedblright\ \textit{IEEE Journal of Selected Areas in
Communications }26 (3), 530-540

\bibitem{NE} Nadakuditi, R.R. and A. Edelman (2008) \textquotedblleft Sample
Eigenvalue Based Detection of High-Dimensional Signals in White Noise Using
Relatively Few Samples\textquotedblright , \textit{IEEE Transactions on
Signal Processing} 56, 2625-2638

\bibitem{NS} Nadakuditi, R.R. and J.W. Silverstein (2010) \textquotedblleft
Fundamental Limit of Sample Generalized Eigenvalue Based Detection of
Signals in Noise Using Relatively Few Signal-Bearing and Noise-Only
Samples\textquotedblright , \textit{IEEE Journal of Selected Topics in
Signal Processing} 4, 468-480.

\bibitem{O1/} Onatski, A. (2009) \textquotedblleft Testing Hypotheses About
the Number of Factors in Large Factor Models\textquotedblright , \textit{%
Econometrica} 77, 1447-1479.

\bibitem{OMH/} Onatski, A., Moreira, M. J., and Hallin, M. (2012)
\textquotedblleft Asymptotic Power of Sphericity Tests for High-dimensional
Data\textquotedblright , manuscript, Economics Faculty, University of
Cambridge.

\bibitem{PPR} Patterson, N., A. L. Price, and D. Reich (2006)
\textquotedblleft Population Structure and Eigenanalysis\textquotedblright , 
\textit{PLoS Genetics} 2 (12), 2074-2093

\bibitem{Pech} P\'{e}ch\'{e}, S. (2003). \textquotedblleft Universality of
local eigenvalue statistics for random sample covariance
matrices,\textquotedblright\ Ph.D. thesis, Ecole Polytechnique F\'{e}d\'{e}%
rale de Lausanne.

\bibitem{PW} Perry, P.O. and P.J. Wolfe (2010) \textquotedblleft Minimax
Rank Estimation for Subspace Tracking\textquotedblright , \textit{IEEE
Journal of Selected Topics in Signal Processing} 4, 504-513

\bibitem{RV/} Ratnarajah, T. and R. Vaillancourt (2005) \textquotedblleft
Complex Singular Wishart Matrices and Applications\textquotedblright , 
\textit{Computers \& Mathematics with Applications} 50, 399--411.

\bibitem{SS} Schreier, P.J., and L. L. Scharf (2010), Statistical Signal
Processing of Complex-Valued Data: The Theory of Improper and Noncircular
Signals, Cambridge University Press.

\bibitem{T} Telatar, E. (1999), \textquotedblleft Capacity of Multi-antenna
Gaussian Channels\textquotedblright , \textit{European Transactions on
Telecommunications} 10 (6), 585-595.

\bibitem{TW} Tracy, C.A., and Widom, H. (2009) \textquotedblleft The
Distributions of Random Matrix Theory and their
Applications\textquotedblright , in V. Sidoravi\v{c}ius (ed.), New Trends in
Mathematical Physics, Springer Science + Business Media B.V., 753-765

\bibitem{TV/} Tulino, A. and Verd\'{u}, S. (2004). Random Matrix Theory and
Wireless Communications. Foundations and Trends in Communications and
Information Theory 1. Now Publishers, Hanover, MA.

\bibitem{vdV/} van der Vaart, A.W. (1998) Asymptotic Statistics, Cambridge
University Press.

\bibitem{W/} Wang, D. (2010) \textquotedblleft The largest eigenvalue of
real symmetric, Hermitian and Hermitian self-dual random matrix models with
rank one external source, part I.\textquotedblright\ arXiv:1012.4144

\bibitem{ZZ} Zinn-Justin, P. and J.-B. Zuber (2003) \textquotedblleft On
some integrals over the U(N) unitary group and their large N
limit\textquotedblright , \textit{Journal of Physics A. }36 (12), 3173-3193
\end{thebibliography}
\end{document}